\documentclass{amsproc}

\newtheorem{theorem}{Theorem}[section]
\newtheorem{lemma}[theorem]{Lemma}

\theoremstyle{definition}

\newtheorem{example}[theorem]{Example}

\theoremstyle{remark}

\numberwithin{equation}{section}

\begin{document}

\title{LECTURES ON SPECIAL LAGRANGIAN SUBMANIFOLDS}

\author{Nigel Hitchin}
\address{Mathematical Institute,
 University of Oxford,
 24-29 St Giles,
 Oxford OX1 3LB,
 England}
\email{hitchin@maths.ox.ac.uk}
\thanks{These notes are derived from lectures given both at the Winter School on Mirror Symmetry in Harvard, January 1999 and the  School on Differential Geometry held at ICTP, Trieste in April 1999. They were initially prepared for ICTP and the author wishes to thank ICTP for permission to publish in these Proceedings.}

\date{30th April 1999}

\begin{abstract}
These notes consist of a  study of special Lagrangian submanifolds of
Calabi-Yau manifolds and their moduli spaces. The particular case of three
dimensions, important in string theory, allows us to introduce the
notion of {\it gerbes}. These offer an appropriate language for describing
many significant features of the Strominger-Yau-Zaslow approach to mirror
symmetry.
\end{abstract}

\maketitle


\newcommand{\oper}[2]{\newcommand{#1}{\mathop{\mathrm{#2}}\nolimits} }
\newcommand{\operlim}[2]{\newcommand{#1}{\mathop{\mathrm{#2}}\limits} }


\newcommand{\lie}[1]{\mathfrak{#1}}
\newcommand{\R}{\mathbf{R}}
\newcommand{\C}{\mathbf{C}}
\newcommand{\Z}{\mathbf{Z}}
\newcommand{\RP}{\R P}
\newcommand{\CP}{\C P}
\newcommand{\Diff}{\operatorname{Diff}}
\newcommand{\Map}{\operatorname{Map}}
\newcommand{\su}{\operatorname{SU}}
\oper{\tr}{tr}
\oper{\ad}{ad}
\oper{\End}{End}
\oper{\Hom}{Hom}
\oper{\sgn}{sgn}
\operlim{\Res}{Res}
\oper{\spinc}{Spin^c}
\oper{\spin}{Spin}
\oper{\so}{SO}
\oper{\gl}{GL}

\section{GERBES}

\subsection{What is a gerbe?}

The word {\it gerbe} is an odd one for English speakers -- it is not a
furry little animal, but derives instead from a French word in more common
use (look at Renoir's painting ``Petite fille \`a la gerbe''). However, it
has in fact existed in the English language for some time. Here's what the
Oxford English Dictionary gives:
\vskip .2cm
\begin{quotation}
{\bf Gerbe.} 1698 [ - Fr. {\it gerbe} wheat-sheaf ]  1. A wheat-sheaf 1808.
2. Something resembling a sheaf of wheat: esp. a kind of firework.
\end{quotation}
\vskip .2cm
\noindent and ``something resembling a sheaf'' is quite close to the mathematical
meaning of the word.

Before giving a definition, it's worthwhile to recognize when we,
as mathematicians, might be in a situation where the language of gerbes
could be relevant. We are basically in gerbe territory (for smooth
manifolds) if any one of the following is being considered:
\vskip .25cm
\begin{itemize}
\item
a cohomology class in $H^3(X,\Z)$
\item
a codimension 3 submanifold $M^{n-3}\subset X^n$
\item
a \v Cech cocycle $[g_{\alpha\beta\gamma}]\in H^2(X, C^{\infty}(S^1))$
\end{itemize}
\vskip .25cm
In the last case, this means a 2-cocycle for the sheaf of germs of
$C^{\infty}$ functions with values in the circle.

The main driving force for the study of gerbes presented here is in trying
to understand the geometry of mirror symmetry for 3-dimensional Calabi-Yau
manifolds. Gerbes are relevant because there is good reason (the study of supersymmetric cycles) to focus
attention on  special Lagrangian submanifolds of 3-dimensional Calabi-Yau
manifolds, and codimension three questions arise in two ways.  Firstly,
such a submanifold is 3-dimensional, and points are of codimension 3; and
secondly the submanifolds themselves are of codimension 3 in the
Calabi-Yau. We shall see both aspects emerging in Lecture 3.

Gerbes were essentially invented by Giraud \cite{Gir} and discussed at
length by Brylinski \cite{Br}. They have also occurred in the work of
Freed and of  Murray.  I have learned much about gerbes from talking to Dan
Freed, and I also  gave my student David Chatterjee the task of trying to
approach gerbes from a concrete point of view. Some of the material of
these lectures  is taken from his thesis \cite{Chat}. The reader should
also look at Brylinski's book.
\vskip .25cm
To understand gerbes, we need to consider the other creatures  in a
hierarchy to which gerbes belong, and here the lowest form of life consists
of circle-valued functions $f:X\rightarrow S^1$. Consider the following
features of such a function:
\vskip .25cm
\begin{itemize}
\item
 a cohomology class in $H^1(X,\Z)$
\item
a codimension 1 submanifold $M^{n-1}\subset X^n$
\item
a \v Cech cocycle $[g_{\alpha}]\in H^0(X, C^{\infty}(S^1))$
\end{itemize}
\vskip .25cm
The cohomology class  is just the pull-back $f^*(x)$ of the generator
$x\in H^1(S^1,\Z)\cong \Z$. If we take the inverse image $f^{-1}(c)$ of a
regular value $c\in S^1$ then this is a codimension one submanifold of $X$
whose homology class is dual to the class $f^*(x)\in H^1(X,\Z)$. Finally,
given an open covering $\{U_{\alpha}\}$ of $X$, a global function
$f:X\rightarrow S^1$ is built up out of local functions
$$g_{\alpha}:U_{\alpha}\rightarrow S^1$$ which satisfy the cocycle
condition $$g_{\beta} g^{-1}_{\alpha}=1\,\,{\rm on}\,\, U_{\alpha}\cap
U_{\beta}$$
Thus the three aspects of gerbes that we identified above all occur here,
but in  degree 1 rather than 3.
\vskip .25cm
The next stage up in the hierarchy consists of a unitary line bundle $L$,
or its principal $S^1$-bundle of unitary frames. Here we have:
\vskip .25cm
\begin{itemize}
\item
a cohomology class in $H^2(X,\Z)$
\item
a codimension 2 submanifold $M^{n-2}\subset X^n$
\item
a \v Cech cocycle $[g_{\alpha\beta}]\in H^1(X, C^{\infty}(S^1))$
\end{itemize}
\vskip .25cm
The degree 2 cohomology class is the  first Chern class $c_1(L)$. If we
take a smooth section of $L$ with nondegenerate zeros, then this vanishes
on $X$ at a codimension two submanifold $M$ Poincar\'e dual to $c_1(L)$.
And this time, if we take an open covering $\{U_{\alpha}\}$ over each set
of which $L$ is trivial, we have transition functions
$$g_{\alpha\beta}:U_{\alpha}\cap U_{\beta}\rightarrow S^1$$
which satisfy $g_{\beta\alpha}=g^{-1}_{\alpha\beta}$ and the cocycle condition
$$g_{\alpha\beta}g_{\beta\gamma}g_{\gamma\alpha}=1\,\,{\rm
on}\,\,U_\alpha\cap U_\beta\cap U_\gamma$$
\noindent Now take an open covering of $X$ and a map
 $$g_{\alpha\beta\gamma}:U_\alpha\cap U_\beta\cap U_\gamma\rightarrow S^1$$
 on each threefold intersection with
$$g_{\alpha\beta\gamma}=g^{-1}_{\beta\alpha\gamma}=g^{-1}_{\alpha\gamma\beta}=g^
{-1}_{\gamma\beta\alpha}$$
 and satisfying the cocycle condition
 $$\delta
g=g_{\beta\gamma\delta}g^{-1}_{\alpha\gamma\delta}g_{\alpha\beta\delta}
  g^{-1}_{\alpha\beta\gamma}=1
 \quad  {\rm on} \quad  U_\alpha\cap U_\beta\cap U_\gamma \cap U_\delta$$
 We shall say that this data defines a {\it gerbe}. By this we mean that it
suffices to define a gerbe in the same way that a collection of transition
functions defines a line bundle or a collection of coordinate charts
defines a manifold. It is imperfect only insofar as it doesn't address the
question of the dependence of the definition on  choices, but I don't want
to complicate things by  introducing sheaves of categories as in \cite{Br}.
This is a working definition, and we are  going to make gerbes work for us.
To take an analogy, any introductory course on manifolds begins with the
existence of a specific collection of charts and then goes on to discuss a
maximal atlas of equivalent charts wherein is supposed to reside the
definition of a manifold. On the other hand, the current version of {\it
MathSciNet} lists 60,509 references to ``manifold'' and only 15 of these
contain a reference  to ``maximal atlas". For working purposes a single
collection of charts will often be enough. Note that a gerbe is 
not a manifold, unlike the total space of a principal bundle.  This is a 
necessary consequence of its definition in terms of threefold intersections
of open sets.  Manifolds by definition involve comparing charts only on
twofold intersections.
\vskip .25cm
An obvious feature of  gerbes is that we can take tensor
products of them (products of the cocycles) and pull them back by maps.
Also, it is clear that a gerbe defined this way  has a characteristic class
in $H^3(X,\Z)$. This follows from the long exact cohomology sequence of the
exact sequence of sheaves
$$0\rightarrow \Z \rightarrow C^{\infty}(\R)\stackrel{\exp 2\pi i
x}{\rightarrow}  C^{\infty}(S^1)\rightarrow 1$$
Since $C^{\infty}(\R)$ is a fine sheaf, we have
$$\rightarrow 0\rightarrow H^2(X,C^{\infty}(S^1))\cong H^3(X,\Z)\rightarrow 0$$
We might say that topologically a gerbe is classified by its characteristic
class just as a line bundle is determined topologically by its Chern class.

\subsection{Trivializations of gerbes}

We can't point to a gerbe as a space in the same way that we can for a line
bundle, but there is a fundamental notion which helps to understand them --
a {\it trivialization} of a gerbe.

Recall that for line bundles a trivialization is a non-vanishing section
$s$ of $L$. A unitary trivialization  ($s/\Vert s \Vert$) is also a section
of the principal $S^1$-bundle. In terms of transition functions this is a
collection of circle-valued functions $f_{\alpha}:U_{\alpha}\rightarrow
S^1$ such that, on $U_{\alpha}\cap U_{\beta}$,
$$f_{\alpha}=g_{\alpha\beta}f_{\beta}$$
Given another trivialization  $f'_{\alpha}:U_{\alpha}\rightarrow S^1$ we
have  $f'_{\alpha}=g_{\alpha\beta} f'_{\beta}$ and so we find that
$$f'_{\alpha}/f_{\alpha}= f'_{\beta}/f_{\beta}$$
Thus $f'_{\alpha}/f_{\alpha}$ is the restriction of a global function to
$U_{\alpha}$ and we see that the difference of two unitary trivializations
of a line bundle is a global {\it function} $f:X\rightarrow S^1$.

We say now that a trivialization of a gerbe is defined by functions
$$f_{\alpha\beta}=f_{\beta\alpha}^{-1}:U_{\alpha}\cap U_{\beta}\rightarrow
S^1$$
on twofold intersections such that
$$g_{\alpha\beta\gamma}=f_{\alpha\beta}f_{\beta\gamma}f_{\gamma\alpha}$$
In other words a trivialization is a specific representation of the cocycle
$g_{\alpha\beta\gamma}$ as a \v Cech coboundary. Here the difference
between any two trivializations $f_{\alpha\beta}$ and $f'_{\alpha\beta}$ is
given by
$$h_{\alpha\beta}= f'_{\alpha\beta}/f_{\alpha\beta}$$
so that $h_{\beta\alpha}=h_{\alpha\beta}^{-1}$ and
$$h_{\alpha\beta}h_{\beta\gamma}h_{\gamma\alpha}=1$$
In other words the difference of two trivializations of a gerbe is a {\it
line bundle}.
\vskip .25cm
 Over a particular open set $U_0$ in our definition of a
gerbe there is a trivialization. We just define, for $\beta,\gamma\ne 0$,
$f_{\beta\gamma}=g_{0\beta\gamma}$, so that the cocycle condition $\delta
g=1$ gives
$$g_{\beta\gamma\delta}=f_{\beta\gamma}f_{\gamma\delta}f_{\delta\beta}$$
If we supplement this with $f_{0\beta}=1$, we have a trivialization as
defined above.
In this approach our gerbe has a local trivialization over each
$U_{\alpha}$ and so on the intersection $U_{\alpha}\cap U_{\beta}$, we have
two trivializations which differ by a line bundle $L_{\alpha\beta}$. With
this we can define a gerbe  using twofold rather than threefold
intersections. The only difference from this point of view is that we move one step up the
hierarchy and use line bundles on $U_{\alpha}\cap U_{\beta}$ rather than
functions. A gerbe is defined  using this language by:
\vskip .25cm
\begin{itemize}
\item
a line bundle $L_{\alpha\beta}$ on each $U_\alpha \cap U_\beta$
\item
  an isomorphism
$L_{\alpha\beta}\cong L^{-1}_{\beta\alpha}$
\item
a trivialization $\theta_{\alpha\beta\gamma}$ of
$L_{\alpha\beta}L_{\beta\gamma} L_{\gamma\alpha}$ (tensor product) on
$U_{\alpha}\cap U_{\beta}\cap U_{\gamma}$
\item
the trivialization $\theta_{\alpha\beta\gamma}$ satisfies $\delta \theta=1$
on $U_{\alpha}\cap U_{\beta}\cap U_{\gamma}\cap U_{\delta}$
\end{itemize}
\vskip .25cm
With regard to the last condition, the  point to notice is that $\delta
\theta$ is a section of a tensor product of 12 line bundles corresponding
to the 12 ordered pairs of elements from $\{\alpha,\beta,\gamma,\delta\}$.
Using the isomorphisms $L_{\alpha\beta}\cong L^{-1}_{\beta\alpha}$, the
products cancel in pairs and so this becomes a section of the trivial
bundle where the trivial section $1$ is well-defined.
\vskip .25cm

 It is through this definition that we can consider concrete examples of
gerbes.

\begin{example}
\noindent Let $\spin(n)$ be the connected double covering of $\so(n)$. Then
the Lie group $\spinc(n)$ is defined as
$$\spinc(n)=\spin(n)\times_{\pm 1} S^1$$
It has surjective homomorphisms $p_1,p_2$ to $\so(n)$ and $S^1$ respectively.
Now suppose $P$ is a principal $\so(n)$ bundle over a manifold $X$. We can
ask whether there exists a principal $\spinc(n)$ bundle $\tilde P$
associated to $P$ by the homomorphism $p_1$. If such a lift exists, a
particular choice is called a $\spinc$-structure, and
two such lifts are equivalent if there exists an isomorphism compatible
with the projection $p_1$. This is a common enough question in differential
topology. Usually the answer is that if the second Stiefel-Whitney class
$w_2(P)$ is the mod 2 reduction of an integral class in $H^2(X,\Z)$ then a
$\spinc$ structure exists, but is not necessarily unique. Two such
structures differ by a principal $S^1$ bundle, as consideration of the
transition functions will reveal, taking account of the fact that
$S^1\subset \spinc (n)$ is in the centre.

Now restricted to a suitably small open set $U_{\alpha}$, the bundle $P$
has a $\spinc$ structure. Think of a choice of such a structure as a
trivialization of a gerbe over $U_{\alpha}$. Then the difference of two
trivializations over $U_{\alpha}\cap U_{\beta}$ is a principal $S^1$ bundle
-- a unitary line bundle $L_{\alpha\beta}$. The required conditions on
these line bundles to define a gerbe follow from the cocycle conditions for
the transition functions $g_{\alpha\beta}:U_{\alpha}\cap
U_{\beta}\rightarrow \so(n)$ for the principal bundle $P$. The
characteristic class in $H^3(X,\Z)$ arises from the exact sequence of
constant sheaves
$$0\rightarrow \Z \stackrel{\times 2}\rightarrow \Z \rightarrow
\Z/2\rightarrow 0$$
which gives
$$\dots \rightarrow H^2(X,\Z)\stackrel{\times 2}\rightarrow
H^2(X,\Z)\rightarrow H^2(X,\Z/2)\rightarrow H^3(X,\Z)\rightarrow \dots$$
The image of $w_2(P)$ under the (Bockstein) homomorphism in this sequence
is the obstruction to this class being the reduction of an integral class,
but for us it is interpreted simply as the characteristic class of the
gerbe. Of course this is 2-torsion and so a very special case.
\end{example}

\begin{example}

\noindent Let $X$ be an oriented 3-manifold, and choose a point $p\in X$.
Now cover $X$ with two open sets,
$U_1\cong \R^3$ a coordinate neighbourhood of $p$, and
 $U_0=X\backslash \{p\}$. Then
 $$U_0\cap U_1\cong \R^3\backslash \{0\}\cong S^2\times \R$$
 We define a gerbe ${\mathcal G}_p$ by taking the line bundle $L_{01}$ on
$U_0\cap U_1$ as the pull-back from $S^2$ of the line bundle whose first
Chern class is the generator of $H^2(S^2,\Z)$, that is the principal $S^1$
bundle $S^3\rightarrow S^2$. The choice of orientation on $X$ gives an
orientation on $S^2$ and hence a choice of generator. There are only two
open sets in this definition so nothing more needs to be checked.
In this case, if $X$ is compact, the characteristic class of ${\mathcal G}_p$
is the generator of $H^3(X,\Z)\cong \Z$. The construction is the analogue
of defining a holomorphic line bundle ${\mathcal L}_p$ from a point $p$ in a
Riemann surface.
\end{example}

\begin{example}
\noindent Take $M^{n-3}\subset X^n$ to be an oriented codimension 3
submanifold of a compact oriented manifold $X$. We can generalize the
previous construction to this situation. Take coordinate neighbourhoods
$U_\alpha$ of $X$ along $M$, then
 $U_\alpha\cong (U_\alpha \cap M)\times \R^3$
Now take
$U_0=X\backslash N(M)$, where $N(M)$ is the closure of a small tubular
neighbourhood of $M$, diffeomorphic to the disc bundle in the normal
bundle. Here $$U_0\cap U_\alpha \cong U_\alpha \cap M\times \{x\in \R^3:
\Vert x \Vert >\epsilon\}$$
and as before we
define $L_{\alpha 0}$ as the pull-back by $x\mapsto x/\Vert x\Vert$ of the
line bundle of degree 1 on $S^2$.
The line bundle $L_{\alpha\beta}= L_{\alpha 0}L^{-1}_{\beta 0}$ is defined on
$$U_{\alpha}\cap U_{\beta}\backslash N(M) \cong (U_{\alpha}\cap U_{\beta}
\cap M)\times \{x\in \R^3: \Vert x \Vert >\epsilon\}\simeq S^2$$
But  by construction $c_1(L_{\alpha\beta})=0$ is zero on $S^2$ and so
extends to a trivial bundle on the whole of $U_{\alpha}\cap U_{\beta}$.
This provides us with all our line bundles to define a gerbe. This is the
analogue for gerbes of the construction of a holomorphic line bundle from a
divisor, and the characteristic class of the gerbe in $H^3(X,\Z)$ is just
the class Poincar\'e dual to the homology class of the submanifold
$M^{n-3}$.
\vskip .25cm
\noindent Note that if we restrict this gerbe to the submanifold $M$ we
used to define it, we get the $\spinc(3)$ gerbe of the normal bundle $N$.
This follows because $\spinc(3)\cong$ U(2) (since $\spin(3)\cong \su(2)$).
Viewing $S^2$ as $\CP^1$, then a $\spinc(3)$ structure on a rank 3
orthogonal vector bundle $V$ is the same thing as finding a rank 2 complex
vector bundle $E$ such that the projective bundle $P(E)$ is isomorphic to
the unit sphere bundle $S(V)$. Given $E$, there is a tautological line
bundle over $P(E)$ whose fibre at $x\in P(E)$ consists of the
one-dimensional space  $x$, and its dual has first Chern class 1.
Conversely a choice of line bundle over $S(V)$ gives the vector bundle $E$.
But now a tubular neighbourhood of $M$ is diffeomorphic to the normal
bundle, so the choice of the line bundle $L_{0\alpha}$ in the definition of
the gerbe fixes a line bundle over $S(N)\vert_{U_{\alpha}}$, and this is
a local trivialization of the $\spinc(3)$ gerbe of $N$. The line bundle
$L_{\alpha\beta}$ is the difference of these two on $U_{\alpha}\cap U_{\beta}$
and so we get globally the $\spinc(3)$ gerbe of the normal bundle.
\end{example}

\subsection{Connections}

Gerbes are not just topological objects: we can do differential geometry
with them too. We shall next describe what a connection on a gerbe is. To
begin with, let's look at a connection on a line bundle which is given by
transition functions
$$g_{\alpha\beta}:U_{\alpha}\cap U_{\beta}\rightarrow S^1\subset \C^*$$
A connection consists of 1-forms $A_{\alpha}$ defined on $U_{\alpha}$ such
that on a twofold intersection $U_{\alpha}\cap U_{\beta}$ we have
$$iA_{\beta}-iA_{\alpha}=g^{-1}_{\alpha\beta}dg_{\alpha\beta}$$
Since $d(g_{\alpha\beta}^{-1}dg_{\alpha\beta})=0$,
$dA_{\beta}-dA_{\alpha}=0$ on $U_{\alpha}\cap U_{\beta}$ and so there is a
global closed 2-form $F$, the {\it curvature},  such that
$$F\vert_{U_{\alpha}}=dA_{\alpha}$$
There is a converse to this, which arises in particular in the study of
geometric quantization (see \cite{Wood}). Suppose $F$ is a closed 2-form on
a manifold $X$ and we take a finite covering by open sets such that all of
their intersections are contractible (e.g. convex neighbourhoods using a
Riemannian metric). Then by the Poincar\'e lemma we have 1-forms
$A_{\alpha}$ defined on $U_{\alpha}$ such that
$F\vert_{U_{\alpha}}=dA_{\alpha}$, and since $dA_{\beta}-dA_{\alpha}=F-F=0$
on $U_{\alpha}\cap U_{\beta}$, applying the Poincar\'e lemma again we have
functions $f_{\alpha\beta}:U_{\alpha}\cap U_{\beta}\rightarrow \R$ with
$A_{\beta}-A_{\alpha}=df_{\alpha\beta}$. Finally if the de Rham cohomology
class $[F]/2\pi\in H^2(X,\R)$ is the image of an {\it integral} class then
there is a choice of $f_{\alpha\beta}$ such that
$g_{\alpha\beta}:U_{\alpha}\cap U_{\beta}\rightarrow S^1$
defined by  $g_{\alpha\beta}=\exp (if_{\alpha\beta})$ is a cocycle. Now we
can reinterpret what we have done as saying that $F$ is the curvature of a
unitary connection on a line bundle with transition functions $g_{\alpha\beta}$.
\vskip .25cm
\noindent Repeatedly using the Poincar\'e lemma is just the standard way of relating
de Rham to \v Cech cohomology, and we can do the same for a closed 3-form
$G$:
\begin{eqnarray*}
G\vert_{U_\alpha}&=&dF_\alpha\\
F_\beta-F_\alpha&=&dA_{\alpha\beta}\\
A_{\alpha\beta}+A_{\beta\gamma}+A_{\gamma\alpha}&=&df_{\alpha\beta\gamma}
\end{eqnarray*}
Moreover if $[G]/2\pi\in H^3(X,\R)$ is the image of an integral class we
can again write $g_{\alpha\beta\gamma}=\exp (if_{\alpha\beta\gamma})$ and
obtain a cocycle. This permits us to define a connection on a gerbe defined
by a cocycle $g_{\alpha\beta\gamma}:U_\alpha\cap U_\beta\cap
U_\gamma\rightarrow S^1$ by forms which satisfy
\begin{eqnarray*}
G\vert_{U_\alpha}&=&dF_\alpha\\
F_\beta-F_\alpha&=&dA_{\alpha\beta}\\
iA_{\alpha\beta}+iA_{\beta\gamma}+iA_{\gamma\alpha}&=&g^{-1}_{\alpha\beta\gamma}
dg_{\alpha\beta\gamma}
\end{eqnarray*}
We call the closed 3-form $G$ the {\it curvature} of the gerbe connection.
Not surprisingly, it represents in real cohomology the characteristic class
of the gerbe.
\vskip .25cm
The $A_{\alpha\beta}$ are 1-forms on twofold intersections, and we could if
we wanted interpret them as connection forms for connections on  line
bundles. In fact, let us  adopt the line bundle definition of a gerbe:
\vskip .25cm
\begin{itemize}
\item
a line bundle $L_{\alpha\beta}$ on each $U_\alpha \cap U_\beta$
\item
  an isomorphism
$L_{\alpha\beta}\cong L^{-1}_{\beta\alpha}$
\item
a trivialization $\theta_{\alpha\beta\gamma}$ of
$L_{\alpha\beta}L_{\beta\gamma}L_{\gamma\alpha}$ on
$U_{\alpha}\cap U_{\beta}\cap U_{\gamma}$ such that
\item
$\delta \theta=1$ on $U_{\alpha}\cap U_{\beta}\cap U_{\gamma}\cap U_{\delta}$
\end{itemize}
\vskip .25cm
 We don't need these open sets to be contractible, and we can define a
connection in this formalism as follows:
\vskip .25cm
\begin{itemize}
\item
a connection $\nabla_{\alpha\beta}$ on $L_{\alpha\beta}$ such that
\item
$\nabla_{\alpha\beta\gamma}\theta_{\alpha\beta\gamma}=0$
\item
a 2-form $F_\alpha\in \Omega^2(U_{\alpha})$ such that on $U_{\alpha}\cap
U_{\beta}$, $F_\beta-F_\alpha=F_{\alpha\beta}=$ the curvature of
$\nabla_{\alpha\beta}$
\end{itemize}
\vskip .25cm
Here $\nabla_{\alpha\beta\gamma}$ is the connection on
$L_{\alpha\beta}L_{\beta\gamma}L_{\gamma\alpha}$ induced by the given
connections on the line bundles $L_{\alpha\beta}$.
\vskip .25cm
We shall say that a connection on a gerbe is {\it flat} if its curvature
$G$ vanishes. Under these circumstances $dF_{\alpha}=0$, so using a
contractible covering,  we can write $F_{\alpha}=dB_{\alpha}$ on
$U_{\alpha}$
 and then on $U_{\alpha}\cap U_{\beta}$,
$$F_{\beta}-F_{\alpha}=dA_{\alpha\beta}=d(B_{\beta}-B_{\alpha})$$
This in turn implies
$$ A_{\alpha\beta}-B_{\beta}+B_{\alpha}=df_{\alpha\beta}$$
But from the definition of connection
$$iA_{\alpha\beta}+iA_{\beta\gamma}+iA_{\gamma\alpha}=g_{\alpha\beta\gamma}^{-1}
dg_{\alpha\beta\gamma}$$
and so
$$d(if_{\alpha\beta}+if_{\beta\gamma}+if_{\gamma\alpha}-\log
g_{\alpha\beta\gamma})=0$$
Of course $\log g$ is only defined modulo $2\pi i \Z$ so what we have here
is a collection of constants
$$c_{\alpha\beta\gamma}\in  2\pi {\bf R/Z}$$
The 2-cocycle $c_{\alpha\beta\gamma}/2\pi$ represents a \v Cech class in
 $ H^2(X,{\bf R/Z})$ which we call the {\it holonomy} of the connection.
 \vskip .25cm
 \noindent For purposes of comparison, the holonomy of a flat connection on
a unitary line bundle is given by parallel translation around closed loops.
It defines a homomorphism from $\pi_1(X)$ to $S^1=\R/\Z$, or equivalently
an element of $H^1(X,\R/\Z)$. If we wish, then we can think of a gerbe with
flat connection as having holonomy around a closed surface $S\subset X$:
since $H^2(S,\R/\Z)\cong \R/\Z$, we get an angle of  holonomy around each
surface. In fact, since the curvature is a 3-form, {\it any} connection is
flat on a surface $S$, and so has holonomy.
\vskip .25cm
For a line bundle, when the holonomy is trivial we get a covariant constant
trivialization of the bundle. What happens for a gerbe?
If the holonomy is trivial, then $c_{\alpha\beta\gamma}$ is a coboundary,
so that there are constants $k_{\alpha\beta}\in 2\pi\R/\Z$ such that
$$c_{\alpha\beta\gamma}=k_{\alpha\beta}+k_{\beta\gamma}+k_{\gamma\alpha}$$
This means that if we put
$$\exp(if_{\alpha\beta}-ik_{\alpha\beta})=h_{\alpha\beta}$$
then
$$h_{\alpha\beta}h_{\beta\gamma}h_{\gamma\alpha}=g_{\alpha\beta\gamma}$$
and so we have a trivialization of the gerbe, which we can call a covariant
constant trivialization or, more briefly, a {\it flat trivialization}. This
is not unique. Suppose we have a second flat trivialization
 $h'_{\alpha\beta}$ then
$g_{\alpha\beta}=h'_{\alpha\beta}/h_{\alpha\beta}$
defines a line bundle $L$. Now
\begin{eqnarray*}
iB_{\beta}-iB_{\alpha}-iA_{\alpha\beta}&=& d\log h_{\alpha\beta}\\
iB'_{\beta}-iB'_{\alpha}-iA_{\alpha\beta}&=& d\log h'_{\alpha\beta}
\end{eqnarray*}
so
$$i(B'-B)_{\beta}-i(B'-B)_{\alpha}=d\log g_{\alpha\beta}$$
and this defines  a  connection on $L$. By definition of $B_{\alpha}$ and
$B'_{\alpha}$,
$$F_{\alpha}=dB_{\alpha}=dB'_{\alpha}$$
so the curvature of this connection is $d(B'_{\alpha}-B_{\alpha})=0$. Thus
the difference of two flat trivializations of a gerbe is a flat line
bundle. We should compare this in the hierarchy with the observation that a
covariant constant unitary trivialization of a line bundle is not unique:
any two differ by a constant function to the circle.
\vskip .25cm
We can apply these ideas to a loop in $X$: a smooth map $f:S^1\rightarrow
X$. Since $S^1$ is only 1-dimensional, the pull back of a gerbe with
connection to the circle is flat and has trivial holonomy. Thus we have
flat trivializations. Suppose we identify flat trivializations if they
differ by a flat line bundle with trivial holonomy. Call the space of
equivalence classes the moduli space of flat trivializations. Then to each
loop we have a space which is acted on freely and transitively by the
moduli space of flat line bundles $H^1(S^1,\R/\Z)\cong S^1$. In other words
we have a principal $S^1$ bundle over the loop space $\Map(S^1,X)$.
Moreover, take a path in the loop space
$$F:[0,1]\times S^1\rightarrow X$$
Since $[0,1]\times S^1 $ is 2-dimensional the pulled back gerbe is flat and
since $[0,1]\times S^1\simeq S^1$ its holonomy is trivial and moreover
there is a canonical isomorphism between the moduli space of flat
trivializations of the gerbe on $\{0\}\times S^1$ and $\{1\}\times S^1$.
This defines the notion of parallel translation in the principal bundle,
and hence a connection on the line bundle over the loop space.

This can be a useful way of viewing gerbes with connection: the curvature
3-form $G$ gets transgressed to a 2-form on the loop space, the curvature
of the line bundle connection defined above. However, to reformulate gerbes
in the more familiar language of line bundles we have had to go to infinite
dimensional spaces. Moreover, the connection is not arbitrary. If we
decompose a surface as a loop of loops in two different ways, the holonomy
is the same. It is like saying that a unitary line bundle with connection
defines a {\it function} $H:\Map(S^1,X)\rightarrow S^1$, the holonomy of
the connection around the loop. This too is a special function: the value
of the map is 1 on any loop which retraces its steps.

\subsection{Examples of connections}

Take our first example of a gerbe, the $\spinc$ gerbe for a principal
$\so(n)$ bundle $P$. The group $\spinc(n)=\spin(n)\times_{\pm 1} S^1$ has a
homomorphism $p_2:\spinc(n)\rightarrow S^1$ which means that any choice of
$\spinc$ structure over $U_{\alpha}$ defines a unitary line bundle
$L_{\alpha}$. The line bundle $L_{\alpha\beta}$ on $U_{\alpha}\cap
U_{\beta}$ then satisfies
$$L^2_{\alpha\beta}=L_{\beta}L_{\alpha}^{-1}$$
From this it is clear that the square of the $\spinc$ gerbe is trivial because we can change the trivialization on $U_{\alpha}$ by the line bundle $L_{\alpha}$ to make $L^2_{\alpha\beta}$ the trivial bundle.
Its characteristic class is therefore a 2-torsion class. If we choose a connection on
each $L_{\alpha}$, with curvature $F_{\alpha}$, this defines a connection
on $L^2_{\alpha\beta}$ with curvature $F_{\beta}-F_{\alpha}$. We can  find
connections $\nabla_{\alpha\beta}$ on $L_{\alpha\beta}$ satisfying the
conditions above with curvature $F_{\alpha\beta}$ and
$$F_{\alpha\beta}=\frac{1}{2}(F_{\beta}-F_{\alpha})$$
Since $F_{\alpha}$ is the curvature of a connection, $dF_{\alpha}=0$ and so
the gerbe connection is flat. Its holonomy is the image of $w_2(P)$ under
the natural homomorphism
$$H^2(X, \Z/2)\rightarrow H^2(X,\R/\Z)$$
\vskip .25cm
\noindent Now consider our second example, the gerbe ${\mathcal G}_p$ for a
point $p\in X$ where $X$ is a compact 3-manifold. Choose a Riemannian
metric on $X$ with volume form $V$ and  total volume $2\pi$.
\vskip .25cm
We shall use the language of currents, or distributional forms. On a
manifold $X^n$ a smooth $p$-form $\alpha$ defines a linear form on
$\Omega^{n-p}(X)$ by
$$\langle \alpha, \beta\rangle=\int_X\alpha \wedge \beta$$
and then by Stokes' theorem
$$\langle d\alpha, \beta\rangle=(-1)^{p+1}\langle \alpha, d\beta\rangle$$
The theory of currents extends this to more general linear forms on
$\Omega^{n-p}(X)$. In particular if $f:Y\rightarrow X$ is a smooth map of
an $(n-p)$-manifold, possibly with boundary, then
$$\beta \mapsto \int_Yf^*\beta$$
is a current. An example is the case $p=n$, and the Dirac delta function
$\delta_x$ of a point $x\in X$. There is a de Rham  and Hodge theory of
currents (see \cite{dR}).
\vskip .25cm
 Now the homology class of a point is dual to the de Rham cohomology class
$V/2\pi$, and so in the language of currents, for our point $p\in X^3$, the
form $V-2\pi \delta_p$ is cohomologous to zero.
It follows from the Hodge theory of currents that we can find a current $H$
such that
$$\Delta H=V-2\pi\delta_p$$
The current is unique modulo a global harmonic 3-form (i.e. a constant
multiple of $V$) and by elliptic regularity  $H$ is a 3-form which is
smooth except at $p$ where the function $\phi=\ast H$ has a singularity of
the form
$$\phi=-\frac{1}{2 r}+\dots$$
Now take a local smooth 3-form $H_1$ on the coordinate neighbourhood $U_1$
such that
$$\Delta H_1=V$$
and define the 2-forms $F_0=d^*H$ on $U_0=X\backslash\{p\}$ and
$F_1=d^*H_1$ on $U_1$. Note that $F_0$ is now independent of the choice of
$H$. We have
$dF_0=dd^*H=\Delta H=V$ on $U_0$ and $dF_1=dd^*H_1=V$ on $U_1$. To produce
a connection with curvature $V$ we need to identify $(F_1-F_0)$ as the curvature of a
connection on a line bundle on $U_0\cap U_1\cong \R^3\backslash \{0\}$ with
first Chern class the generator of $H^2(\R^3\backslash \{0\})\cong \Z$. So
take a closed ball $B$ centred on $p$ in the coordinate neighbourhood $U_1$
and $\varphi$ a smooth function of compact support in $U_1$ which is
identically 1 in a neighbourhood of $B$. By the definition of $H$ and $H_1$,
$$\langle dd^*H, \varphi\rangle=\int_{U_1}\varphi V-2\pi\varphi(p)$$
$$\int_{U_1} \varphi dd^*H_1=\int_{U_1}\varphi V$$
so subtracting,
$$\langle d(d^*H-d^*H_1),\varphi\rangle =-2\pi$$
since $\varphi(p)=1$. But $d^*H-d^*H_1$ is a smooth {\it closed} 2-form
outside $p$, so this is the same as restricting to $B$, and there
$\varphi=1$, so we have
$$-2\pi=\langle B, d(d^*H-d^*H_1)\rangle=\langle
dB,(d^*H-d^*H_1)\rangle=\int_{\partial B}(d^*H-d^*H_1)$$
Thus  $(F_1-F_0)/2\pi=d^*(H_1-H)/2\pi$ has integral cohomology class, a
generator of the cohomology group, and so $(F_1-F_0)$ is the curvature of
the required connection on the line bundle $L_{01}$. We therefore have a
gerbe connection on ${\mathcal G}_p$ with curvature $V$.

\subsection{An infinitesimal viewpoint}

Although I have been emphasizing a down to earth point of view in this
discussion on gerbes, we finish here with an attempt to rid the notion of a
connection on a gerbe from any dependence on a covering. The reader might
wish to skip this and go on to the applications of gerbes.

From the description of gerbe connections above, it looks as if a
trivialization of a connection defines a connection 2-form $F_{\alpha}$
just as a trivialization of a line bundle defines a connection 1-form
$A_{\alpha}$. This is not quite true. To see what really happens, we
should study connections on the trivial gerbe, that is where
$g_{\alpha\beta\gamma}=1$. In this case, a connection is defined by
\begin{eqnarray}
F_\beta-F_\alpha&=&dA_{\alpha\beta}\label{first}\\
A_{\alpha\beta}+A_{\beta\gamma}+A_{\gamma\alpha}&=&0\label{second}
\end{eqnarray}
The second equation (\ref{second}) says we have a 1-cocycle with values in
$\Omega^1$. It defines a vector bundle $E$ which is an extension of the
trivial bundle by the cotangent bundle:
$$0\rightarrow T^*\rightarrow E\stackrel{\pi}\rightarrow 1\rightarrow 0$$
To be concrete, this is described by transition functions which are matrices
$$
\left(
\begin{array}{ccccc}
1&0&0&\dots&0\\
                 a_1&1&0&\dots&0\\
                 a_2&0&1&\dots&0\\
                \dots&\dots&\dots&\dots&\dots\\
                   a_n&0&0&\dots&1
                   \end{array}
                   \right)
                 $$
 where $A_{\alpha\beta}=a_1dx_1+a_2dx_2+\dots+a_ndx_n$.
  The first  equation (\ref{first}) defines a splitting of the extension
$\Lambda^2T^*\rightarrow E'\rightarrow 1$ determined from the exterior
derivative of $A_{\alpha\beta}$. In concrete terms, we can consider the
extension obtained from $E$ just by tensoring with $T^*$
$$0\rightarrow T^*\otimes T^*\rightarrow T^*\otimes E\rightarrow
T^*\rightarrow 0$$
and if we define $T^*\wedge E$ by $T^*\otimes E/Sym^2 T^*$ then we have
$$0\rightarrow\Lambda^2 T^*\rightarrow T^*\wedge E\stackrel{\pi}\rightarrow
T^*\rightarrow 0$$
and the  splitting of the exterior derivative extension determined by
$F_{\alpha}$ can be interpreted as a first order linear differential
operator
$$D:C^{\infty}(E)\rightarrow  C^{\infty}(T^*\wedge E)$$
with the properties that for a section $e\in C^{\infty}(E)$
\vskip .25cm
\begin{itemize}
\item
$D(fe)=df\wedge e+fDe$
\item
if $e\in C^{\infty}(T^*)\subset C^{\infty}(E)$ then $De = de\in
C^{\infty}(\Lambda^2T^*)\subset C^{\infty}(T^*\wedge E)$
\item
 if $\pi(e)=f$ then $\pi(De)=df$
 \end{itemize}
 \vskip .25cm
 Any two such operators on the same bundle $E$ differ by a zero order term:
 $$(D_1-D_2)(e)=F\pi(e)$$
 where $F$ is a 2-form on $X$.
 \vskip .25cm
\noindent There is a natural class of  examples of such extensions  and
differential operators.  Let $L$ be a line bundle and $J_1(L)$ the bundle
of 1-jets of sections: its fibre at a point $x\in X$ is $C^{\infty}(L)/K_x$
where $K_x$ is the vector space of   sections which vanish together with
their first derivative at $x$.  Given any section $s \in C^{\infty}(L)$ its
image in  $C^{\infty}(L)/K_x$ as $x$ varies defines a section $j_1(s)$ of
$J_1(L)$ -- the 1-jet of $s$. Evaluating the section at $x$ gives a
homomorphism $\pi:J_1(L)\rightarrow L$ with kernel $L \otimes T^*$. Thus
$J_1(L)$ is an extension
$$ 0\rightarrow L\otimes T^*\rightarrow     J_1(L)\rightarrow L\rightarrow 0$$
and
$L^{-1}J_1(L)$ is an extension:
\begin{equation}
 0\rightarrow T^*\rightarrow     L^{-1}J_1(L)\rightarrow 1\rightarrow 0
 \label{ext}
 \end{equation}
If $L$ is a complex line bundle then this is strictly speaking complex,
with the extension class defined by
$A_{\alpha\beta}=g_{\alpha\beta}^{-1}dg_{\alpha\beta}$ but a unitary
structure expresses it as a real bundle since then $A_{\alpha\beta}$ is
imaginary.  Another way of looking at it, as in \cite{A}, is to take the
principal circle bundle $P$ of the line bundle $L$. Then $T^*P/S^1\cong
L^{-1}J_1(L)$. From either point of view a  unitary connection on $L$ is
just a splitting of the extension (\ref{ext}).

Now the bundle $L^{-1}J_1(L)$ has a distinguished family of local sections,
namely those of the form $s^{-1}j_1(s)$ for a non-vanishing local
section $s$ and there is a differential operator $D_{L}$ as above such that
any local solution to $D_{L}(e)=0$ with $\pi(e)=1$ is of the form
$e=s^{-1}j_1(s)$.
\vskip .25cm
We see that a connection on the trivial gerbe is defined by a pair
$(E,D)$. A choice of splitting of the extension $E\cong T^*\oplus 1$ allows
us to define $D$ by a 2-form $F$
$$D(\alpha, f)=(d\alpha+Ff, df)$$
The curvature of this gerbe connection is $dF$, but now  the connection
``form'' is $F$ only after choosing a splitting.
\vskip .25cm
If we have a trivialization of the trivial gerbe, this is given by
$h_{\alpha\beta}$ satisfying
$h_{\alpha\beta}h_{\beta\gamma}h_{\gamma\alpha}=1$
and
$$iA'_{\alpha\beta}=iA_{\alpha\beta}+h^{-1}_{\alpha\beta}dh_{\alpha\beta},\qquad
F'_{\alpha}=F_{\alpha}$$
can be thought of as defining the same connection relative to the new
trivialization. If we now regard the difference of two trivializations of a
gerbe as a line bundle $L$, then we can say that the bundle of the
connection  $E'$ for one trivialization differs (as an extension) from the
bundle $E$ of another by $L^{-1}J_1(L)$, with the differential operators
coinciding. Thus a connection on a gerbe, relative to local trivializations
over open sets $U_{\alpha}$, is given by pairs $(E_\alpha,D_\alpha)$ such
that on $U_{\alpha}\cap U_{\beta}$
$$E_{\beta}-E_{\alpha}=L_{\alpha\beta}^{-1}J_1(L_{\alpha\beta}),\qquad
D_\beta-D_\alpha=D_{L_{\alpha\beta}}$$
This makes a gerbe connection more like that of a connection on  a line
bundle. The description above of a connection in terms of 2-forms and
connections on line bundles just involves the extra data of splitting all
the extensions. We can push the language further by thinking of the vector
bundle $E$ over a point $x$ as the 1-jet of a line bundle and $E$ together
with its differential operator $D$ as the 2-jet of a line bundle. We know
that the difference of two trivializations is a line bundle, so we would
like to formulate the notion of a connection on a gerbe in the following
infinitesimal form:
\vskip .25cm
\noindent {\it A gerbe connection  extends each trivialization at a point
to a 2-jet of a trivialization. }
\vskip .25cm
\noindent For most purposes this may not be a useful working definition,
but it is certainly compatible with the view of parallel translation in
line bundles as a distinguished way of extending a non-zero section of a
line bundle at a point to the 1-jet of a section.

\section{ MODULI SPACES OF SPECIAL LAGRANGIAN SUBMANIFOLDS}

\subsection {Special Lagrangian submanifolds}

A {\it Calabi-Yau manifold} is a  K\"ahler manifold of complex dimension
$n$ with a covariant constant holomorphic $n$-form. Equivalently it is a
Riemannian  manifold with holonomy contained in $\su(n)$.

All Ricci-flat manifolds with special holonomy share the feature that their
geometry is determined by a collection of closed forms. In this case
we have the K\"ahler 2-form $\omega$ and the real and imaginary parts
$\Omega_1$ and $\Omega_2$ of the covariant constant $n$-form (note we have
here made a choice for $\Omega_1+i\Omega_2$ -- any other will differ by
multiplying by a  complex number). These forms satisfy some  identities:
\vskip .25cm
\begin{itemize}
\item
 $\omega$ is non-degenerate
\item
$\Omega_1+i\Omega_2$ is locally decomposable and non-vanishing
\item
 $\Omega_1 \wedge \omega=\Omega_2 \wedge \omega =0$
\item
 $(\Omega_1+i \Omega_2)\wedge(\Omega_1-i \Omega_2)=c\omega^n$ for a
constant $c$
\item
 $d\omega=0,\quad d\Omega_1=0,\quad d\Omega_2=0$
\end{itemize}
\vskip .15cm
and one can show in fact (see \cite{H}) that these suffice to define the
metric and complex structure.

A submanifold $M$ of a symplectic manifold $Z$ is  Lagrangian if $\omega$
restricts to zero on $M$ and $\dim Z=2 \dim M$. A submanifold of a Calabi-Yau
manifold is  {\it special Lagrangian} if in addition $\Omega_2$ restricts
to zero on $L$. In this case (when we choose $c$ appropriately), the real
part $\Omega_1$ of the holomorphic $n$-form when restricted to $M$ is  the
volume form $V$ of the induced metric on $M$.

Concrete examples of special Lagrangian submanifolds are difficult to find,
and so far consist of  three types:
\vskip .25cm
\begin{itemize}
\item
complex Lagrangian submanifolds of hyperk\"ahler manifolds
\item
fixed points of a real structure on a Calabi-Yau manifold
\item
explicit examples for non-compact Calabi-Yau manifolds
\end{itemize}
\vskip .25cm
In this lecture we shall be looking at 3-dimensional examples, in which
case the hyperk\"ahler case is irrelevant.
If $Z$ is a Calabi-Yau manifold with a real structure --- an antiholomorphic
involution $\sigma$
 for which $\sigma^*\omega=-\omega$ and $\sigma^*\Omega_2=-\Omega_2$ --- then
 the fixed point set (the set of real points of $Z$) is easily seen to be a
 special  Lagrangian  submanifold. Real forms of concrete Calabi-Yau's
occur for example by taking a quintic hypersurface in $\CP^4$ given by a
polynomial with real coefficients.
 \vskip .15cm
 In the non-compact case, Stenzel \cite{Sten} has some concrete
 examples. In particular $T^*S^n$ (with the complex structure of an affine
 quadric $\sum_0^n z_i^2=1$) has a complete Calabi-Yau metric which one can
write down. In this case  the zero section is
 special Lagrangian. This is the sphere $\sum_0^n x_i^2=1$, the fixed point
set of the real structure $z_i\mapsto \bar z_i$. The noncompact fibres are
also special Lagrangian.
\vskip .25cm
Non-explicit examples can be shown to exist by deformation. The fundamental
tool here is the theorem of McLean \cite{Mac}. This shows that given one
compact special Lagrangian submanifold $M_0$, there is a local finite
dimensional moduli space  of deformations whose dimension is equal to the
first Betti number $b_1(M_0)$. Thus in Stenzel's example the zero section
of $T^*S^n$, being simply connected, is rigid. Nevertheless, starting with
a set of real points in a suitable compact Calabi-Yau, and deforming, one
can certainly assert the  existence of compact special Lagrangian
submanifolds. One such approach can be found in \cite{B}. The problem is to
investigate the global structure of the moduli space and so far the
appropriate tools are not available. Despite the lack of explicit examples,
we can still study the differential geometry of the situation.

\subsection {Moment maps and gerbes}

In this lecture I want to give a description of the natural local
differential geometry of the moduli space from the point of view of
symplectic geometry and moment maps.
Gerbes will appear in a natural way here when the complex dimension of the
Calabi-Yau manifold (and hence the real dimension of a special Lagrangian
submanifold) is 3. The original idea is due to S.K. Donaldson \cite{Don}, as  an
approach to one of many problems in differential geometry. We
begin with the observation that the volume form of a special Lagrangian
submanifold $M$ is the restriction of the closed form $\Omega_1$. This
means the total volume of  $M_0$  just entails evaluating the cohomology
class $[\Omega_1]\in H^3(Z, \R)$ on the fundamental class $[M_0]\in
H_3(M_0,\Z)$. Thus when we deform $M_0$ using McLean's theorem, the total
volume remains the same. Now by a theorem of Moser \cite{Mo}, this means
that each submanifold $M$ in the family of deformations is diffeomorphic as
a manifold with volume form to $M_0$.

We start then with a fixed compact oriented 3-manifold $M$ with volume form
$V$, and consider the infinite dimensional space
$$\Map(M,Z)$$
of smooth maps to a Calabi-Yau 3-fold $Z$. A tangent vector to this space
at a map $f:M\rightarrow Z$ is a section $X$ of the pulled back tangent
bundle $f^*(TZ)$. We now define a 2-form on the space of maps by
$$\varphi(X,Y)=\int_M \omega(X,Y)V$$
This is formally  non-degenerate, but also,
\begin{lemma} The 2-form $\varphi$ is closed.
\end{lemma}
\begin{proof}
Consider the evaluation map
$$E:\Map(M,Z)\times M\rightarrow Z$$
defined by $E(f,m)=f(m)$. Let $p$ be the projection from the product to the
first factor,
then the definition of $\varphi$ above is equivalent to
$$\varphi=p_*(E^*\omega \wedge V)$$
where $p_*$ denotes integrating over the fibres $M$. Since $d\omega=0$ and
$dV=0$, and $M$ is compact without boundary,
$$d\varphi=dp_*(E^*\omega \wedge V)=p_*(d(E^*\omega \wedge V))=0$$
\end{proof}
We can thus think of $\Map(M,Z)$ as an infinite-dimensional symplectic
manifold. Moreover, the group $\Diff_V(M)$ of volume-preserving
diffeomorphisms acts on the space $\Map(M,Z)$ preserving the symplectic form. We can
ask whether there is a moment map for this action.
\vskip .25cm
Recall that if a Lie group $G$ acts on a  manifold $N$, there is a
homomorphism $\xi\mapsto X(\xi)$ from the Lie algebra ${\lie g}$ of $G$ to
the Lie algebra of vector fields  and if the group action preserves a
symplectic form $\varphi$ then since $d\varphi=0$
$$0={\mathcal
L}_{X(\xi)}\varphi=d(\iota(X(\xi))\varphi)+\iota(X(\xi))d\varphi=d(\iota(X(\xi))
\varphi)$$
so for each $\xi$, $\iota(X(\xi))\varphi$ is a closed 1-form and hence
locally the differential of a function. If there exists a global
equivariant function
$$\mu:N\rightarrow {\lie g}^*$$
to the dual of the Lie algebra such that
$$d\langle \mu, \xi\rangle=\iota(X(\xi))\varphi$$
then $\mu$ is called a moment map.
\vskip .25cm
To find a moment map for the action above with $N=\Map(M,Z)$ and $G=\Diff_V(M)$, we need to look at the Lie
algebra  ${\lie g}$ of the group $\Diff_V(M)$. This consists of the vector fields $X$
on $M$ such that ${\mathcal L}_XV=0$, or equivalently $d(\iota(X)V)=0$ since
$V$ is of top degree and hence closed. Thus since $V$ is everywhere
non-vanishing, $X\mapsto \iota(X)V$ gives an isomorphism from
the Lie algebra to the space of closed 2-forms. Taking the de
Rham cohomology class we have vector space homomorphisms
\begin{equation}
0\rightarrow \lie{g}_0\rightarrow \lie{g}\rightarrow H^2(M,{\R})\rightarrow 0
\label{seq}
\end{equation}
Now since ${\mathcal L}_XV=0$
$$\iota([X,Y])V=\iota({\mathcal L}_XY)V={\mathcal
L}_X(\iota(Y)V)=d(\iota(X)\iota(Y)V)$$
since $\iota(Y)V$ is closed. Hence the cohomology class of $\iota([X,Y])V$
is trivial, and (\ref{seq}) is a sequence of Lie algebra homomorphisms,
with the trivial Lie bracket on $H^2(M,{\R})$. The Lie algebra ${\lie g_0}$
is the space of exact 2-forms and its dual is
$${\lie g_0}^*=\Omega^1(M)/\ker d$$
We see this through the non-degenerate pairing between 1-forms and 2-forms
$$\langle \alpha, \beta\rangle=\int_M \alpha\wedge \beta$$
so if
$\langle \alpha, \beta\rangle=0$ for all exact $\beta\in \Omega^2(M)$ then
$$0=\int_M \alpha\wedge d\gamma=\int_M d\alpha\wedge \gamma$$
for all $\gamma \in \Omega^1(M)$ which implies $d\alpha=0$. Thus ${\lie
g_0}^*$ is the quotient of all 1-forms modulo closed 1-forms as claimed.
\vskip .25cm
We need to find a Lie group whose Lie algebra is $\lie g_0$, and here
gerbes come to the rescue. Assume that the total volume of $M$ is $2\pi$. Then
$[V/2\pi ]\in H^3(M,\R)$ is an integral class and so there are gerbes with
connection whose curvature is $V$ (call these ``degree 1 gerbes'').
Let ${\mathcal G}$ be a degree 1 gerbe then so is $h^*{\mathcal G}$ if
$h:M\rightarrow M$ is a volume preserving diffeomorphism. Hence $(h^*{\mathcal
G}){\mathcal G}^{-1}$ has curvature zero and so is a flat gerbe. The holonomy
of this gerbe  therefore determines a  map
$$\Diff_V(M)\rightarrow H^2(M,\R/\Z)$$
Let $G\subset \Diff_V(M)$ be the component of the identity, then $G$ acts
trivially on the cohomology and the map above restricted to $G$ is a
homomorphism of groups:
$$1\rightarrow G_0\rightarrow G \rightarrow H^2(M,\R/\Z)\rightarrow 1$$
The group $G_0$ can therefore be interpreted as the subgroup of the
identity component of the group of volume preserving diffeomorphisms which
preserves the equivalence class of each degree one gerbe. We shall show
that $G_0$ has a moment map.
\vskip .25cm
The space $\Map(M,Z)$ has in general many components. In particular, for
each map $f:M\rightarrow Z$ we have a distinguished cohomology class
$$[f^*\omega]\in H^2(M,\R)$$
Since we are looking for Lagrangian submanifolds, those for which
$f^*\omega=0$, we need only restrict our search to the subspace
$\Map_0(M,Z)$ for which this class is zero. Thus for each $f\in
\Map_0(M,Z)$,
$$f^*\omega=d\theta$$
where $\theta \in \Omega^1(M)$ is well-defined modulo closed 1-forms. For
each $f$ we thus have a well defined element
$$[\theta_f]\in \Omega^1(M)/\ker d=\lie{g}_0^*$$
\vskip .25cm
\begin{theorem} (see \cite{Don}) The group $G_0$ acting on $\Map_0(M,Z)$ has moment map
$$\mu(f)=[\theta_f]$$
\end{theorem}
\vskip .25cm
\begin{proof}
We first have to find the map $\xi\mapsto X(\xi)$ arising from the action
of $\Diff_V(M)$. For a map $f:M\rightarrow Z$ this action is composition
with the diffeomorphism $h^{-1}\in \Diff_V(M)$:
$$h(f)=f\circ h^{-1}$$
Under the evaluation map, we have
$$E(h(f),h(x))=h(f)(h(x))=f(h^{-1}h(x))=f(x)=E(f,x)$$
Thus $\xi \in \lie{g}$ acts on $\Map_0(M,Z)$ by the vector field $X(\xi)$ where
\begin{equation}
dE(X(\xi),\xi)=0
\label{de}
\end{equation}
Now write $\omega= d\tilde \theta$ in a neighbourhood of $f(M)\in Z$. We have
$$E^*(\tilde \theta)=(\theta_1,\theta_2)$$
and $\theta_1$ is still $\Diff_V(M)$ invariant.
From (\ref{de}) we obtain
\begin{equation}
\theta_1(X(\xi))+\theta_2(\xi)=0
\label{theta}
\end{equation}
Now by the definition of the symplectic form $\varphi$, we have, near $f$,
$$\varphi=p_*(E^*\omega \wedge V)=p_*(E^*d\tilde \theta \wedge
V)=dp_*(E^*\tilde \theta \wedge V)=d\psi$$
and now
\begin{equation}
\iota(X(\xi))\varphi=\iota(X(\xi))d\psi=-d(\iota(X(\xi))\psi)
\label{one}
\end{equation}
since $\psi$ is $\Diff_V(M)$ invariant. But $\psi=p_*(E^*\tilde\theta
\wedge V)$ so from (\ref{theta})
\begin{equation}
\iota(X(\xi))\psi=\int_M\iota(X(\xi))\theta_1 V=-\int_M(\iota(\xi)\theta_2) V
\label{two}
\end{equation}
 But $\theta_2\wedge V=0$ since $M$ is only 3-dimensional, so
$$0=\iota(\xi)(\theta_2\wedge V)=(\iota(\xi)\theta_2)\wedge
V-\theta_2\wedge \iota(\xi)V$$
Thus from (\ref{one}) and (\ref{two}),
$$\iota(X(\xi))\varphi=d\int_M(\iota(\xi)\theta_2) V=d\int_M \theta_2\wedge
(\iota(\xi)V)$$
Thus for $\xi \in \lie {g_0}$, $\iota(\xi)V$ is exact and the class of
$\theta_2$ (which restricted to $\{f\}\times M\in \Map_0(M,Z)\times M$ is
the form $\theta$ for which $f^*\omega=d\theta$) is the moment map.

\end{proof}
\vskip .25cm
From this proposition, we can see that the zero set of the moment map
consists of maps $f$ for which $f^*\omega=d\theta$ where $\theta$ is
closed, i.e. $f^*\omega=0$. Amongst these are the Lagrangian submanifolds.
We consider next the {\it special} Lagrangian condition.

\subsection {K\"ahler quotients}

We used the symplectic form $\omega$ on $Z$ to put a symplectic structure
on $\Map(M,Z)$. We can use the complex structure $I:TZ\rightarrow TZ$ to
define a complex structure on $\Map(M,Z)$. If $X\in C^{\infty}(f^*TZ)$ is a
section of the pulled back tangent bundle, so is $IX$, and formally
speaking this structure is integrable. Now consider the subspace
$${\mathcal S}=\{ f\in \Map_0(M,Z): f^*\Omega_2=0, f^*\Omega_1=V\}$$

\noindent Note that imposing the condition $f^*\Omega_1=V$ implies that
maps in ${\mathcal S}$ are immersions.

\begin{lemma} The subspace ${\mathcal S}$ is a complex submanifold of  $Map(M,Z)$.
\end{lemma}
\begin{proof}
Differentiating the two conditions $f^*\Omega_2=0$ and $f^*\Omega_1=V$, we
see that
 $X\in C^{\infty}(f^*TZ)$ is tangent to ${\mathcal S}$ if and only if
 \begin{equation}
 d(\iota(X)\Omega_1)=0,\qquad d(\iota(X)\Omega_2)=0
 \label{tan}
 \end{equation}
 when pulled back to $M$. But $\Omega_1+i\Omega_2$ is a 3-form on $Z$ of
type $(3,0)$ so that
\begin{equation}
\iota(IX)(\Omega_1+i\Omega_2)=i\iota(X)(\Omega_1+i\Omega_2)
\label{oo}
\end{equation}
and hence
$$\iota(IX)\Omega_1=-\iota(X)\Omega_2,\qquad
\iota(IX)\Omega_2=\iota(X)\Omega_1$$
From (\ref{tan}) we see that if $X$ is tangent to ${\mathcal S}$ then so is
$IX$, so that ${\mathcal S}$ is complex.
\end{proof}
Now ${\mathcal S}$ is a complex submanifold of a K\"ahler manifold and hence
K\"ahler. Moreover it is invariant under $\Diff_V(M)$, so we can consider
the action of $G_0\subset \Diff_V(M) $ on it. The zero set of the moment
map  for this action on ${\mathcal S}$ consists of the space of immersions of
$M$ into $Z$ such that
$$f^*\omega=0,\quad f^*\Omega_1=V,\quad f^*\Omega_2=0$$
 so that the embeddings in this family are precisely the special Lagrangian
submanifolds which are diffeomorphic as manifolds with volume form to
$(M,V)$.
They are, however, represented as maps, in other words as parametrized
submanifolds. To remove the parametrization would mean taking the quotient
of
${\mathcal S}\cap \mu^{-1}(0)$ by $\Diff_V(M)$.
\vskip .25cm
As is well-known, if a group acts on a symplectic manifold with moment map
$\mu$ and acts freely and discontinuously on $\mu^{-1}(0)$, then the
quotient
$$\mu^{-1}(0)/G$$
is again a symplectic manifold, the symplectic quotient.
Furthermore, taking the symplectic quotient of a K\"ahler manifold by a
group which preserves both the complex structure and K\"ahler form yields
again a K\"ahler manifold, using the quotient metric. Now
$${\mathcal S}\cap \mu^{-1}(0)$$
is the zero set of the moment map for $G_0$ acting on the K\"ahler manifold
${\mathcal S}$, and so ${\mathcal S}\cap \mu^{-1}(0)/G_0$ is, formally speaking, a
manifold with a natural K\"ahler metric.
\vskip .25cm
 What is it? Using gerbes we can say precisely. Since $G_0$ was the
subgroup preserving equivalence classes of degree one gerbes, this quotient
is
 {\it the moduli space ${\mathcal M}$ of pairs $(M, {\mathcal G})$ where $M\subset
Z$ is a special Lagrangian submanifold and ${\mathcal G}$  a degree one gerbe
on $M$.}

\subsection{The K\"ahler metric}

To a certain degree what we did in the previous section was formal. It has
to be supported by a theorem which asserts that there really is a manifold
which is the moduli space so described. But this is where McLean's work
comes in.
\vskip .25cm
 To begin with consider the tangent space of the moduli space  ${\mathcal M}$
as described by the symplectic quotient construction above. Tangent vectors
are sections $X \in C^{\infty}(f^*TZ)$ which satisfy three conditions. The
first expresses the fact that we have a vector tangent to the zero set of
the moment map for $G_0$ and then two more which describe the tangent space
to the K\"ahler submanifold ${\mathcal S}$:
\vskip .25cm
\begin{itemize}
\item
$d(\iota(X)\omega)=0$
\item
$d(\iota(X)\Omega_1)=0$
\item
$d(\iota(X)\Omega_2)=0$
\end{itemize}
\vskip .25cm
All three expressions are to be interpreted as forms restricted to the
special Lagrangian submanifold $f(M)$ (which we shall call $M$). The
tangent space of the moduli space is formally the quotient of this space of
sections by the tangent space to the orbit of $G_0$.

To analyse this closer, note that since $\omega$ vanishes on $M$, then $M$
itself has no complex tangent vectors, so the orthogonal complex structure
$I:TZ\rightarrow TZ$ on the Calabi-Yau manifold $Z$ maps $TM\subset TZ$ to
its orthogonal complement, the normal bundle $N \subset TZ$. Now decompose
each section $X$ of $TZ$ over $M$ into its normal and tangential parts:
$$X=X_{ n}+X_{ t}$$
The special Lagrangian condition says that $\omega$ and $\Omega_2$ restrict
to zero on $M$, so that $\iota(X_{ t})\omega$ and $\iota(X_t)\Omega_2$ are
both zero on $M$. The first and third conditions thus only depend on $X_n$.
We use this to prove the following:

\begin{lemma} With respect to the induced metric on $M$,
$\iota(X)\omega=-\ast \iota(X)\Omega_2$. Since $d(\iota(X)\omega)=0$
 and
$d(\iota(X)\Omega_2)=0$ it follows that $\iota(X)\omega$ is a harmonic
1-form on $M$.
\end{lemma}
\begin{proof} On an oriented Riemannian manifold with volume form $V$, let
$X$ be the vector field dual to the 1-form $\alpha$ (i.e.
$\alpha(Y)=g(X,Y)$ for all vector fields $Y$). Then
$$\ast \alpha=\iota(X)V$$
Now take $\alpha=\iota(X_n)\omega$. The vector field dual to this in the
K\"ahler metric on $Z$ is $IX_n$, which is tangent to $M$, so this is its
dual on $M$ in the induced metric. Now since $\Omega_1$ restricts to the
volume form $V$ on $M$,
$$\iota(IX_n)V=\iota(IX_n)\Omega_1=-\iota(X_n)\Omega_2$$
from (\ref{oo}), which completes the proof.
\end{proof}
\vskip .25cm
 From the lemma, the normal components of $X$ have the property that
$\iota(X_n)\omega$ lies in the finite-dimensional space of harmonic
1-forms. By Hodge theory this has dimension $b_1(M)$, the first Betti
number of the compact manifold $M$. Now consider the second condition
$d(\iota(X)\Omega_1)=0$. Since from (\ref{oo})
$\iota(X_n)\Omega_1=\iota(IX_n)\Omega_2=0$ on $M$, this is a condition only
on the tangential part $X_t$, but since $\Omega_1$ restricts to the volume
form $V$, $d(\iota(X)\Omega_1)=0$  is simply the condition that $X_t$ is
volume preserving. As we know that $\Diff_V(M)$ acts on our space, $X_t$ is
thus an arbitrary volume preserving vector field. Our moduli space is the
quotient by the action of $G_0$, whose Lie algebra is defined by the exact
2-forms, so the quotient consists of closed forms modulo exact forms -- the
cohomology group  $H^2(M, \R)$. We obtain therefore a description of the
tangent space of the moduli space as an extension
$$0\rightarrow H^2(M,\R)\rightarrow T\rightarrow H^1(M,\R)\rightarrow 0$$
In particular since $\dim M=3$, $\dim H^2(M,\R)=\dim H^1(M,\R)$ and so the
supposed tangent space has real dimension $2b_1(M)$, or complex dimension
$b_1(M)$.
\vskip .25cm
Passing from the formal to the actual situation, we use McLean's result,
which shows that every harmonic 1-form arises as a tangent vector for a one
parameter family of deformations of a given special Lagrangian submanifold:

\begin{theorem} \cite{Mac}  A normal vector field $X$ to a compact
special Lagrangian submanifold $M$ is the deformation vector field to a normal
deformation through special Lagrangian submanifolds if and only if the
corresponding 1-form $IX$ on $M$ is harmonic. There are no obstructions to
extending a first order deformation to an actual deformation and the
tangent space to such
deformations can be identified through the cohomology class of the harmonic
form with $H^1(M,{\bf R})$.
\end{theorem}
\vskip .25cm
With this result we have the structure of a $b_1(M)$-dimensional real
manifold on the moduli space of special Lagrangian deformations of $M$. We
know already that the equivalence classes of gerbes form a torus, so we
have a K\"ahler manifold of complex dimension $b_1(M)$ as our symplectic
quotient. Let us look closer at its differential-geometric properties.
\vskip .25cm
First we note that the full group $\Diff_V(M)$ acts on ${\mathcal S}$
preserving the symplectic form and complex structure and hence the metric.
Thus the quotient group
$\Diff_V(M)/G_0$ acts on ${\mathcal M}$ holomorphically and symplectically and
hence also its connected subgroup
$$G/G_0\cong H^2(M,\R/\Z)$$
which is a torus $T^n$ of dimension $n=b_1(M)$. It acts freely on ${\mathcal
M}$ and its quotient is the moduli space ${\mathcal B}$ of special Lagrangian
submanifolds.
\vskip .25cm
The group $\Diff_V(M)$ also acts as reparametrizations of the submanifold
$M$, so that the tangent vectors to its orbits in $\Map(M,Z)$ are sections
$X\in C^{\infty}(f^*TZ)$ which lie in $TX$. From the definition of the
symplectic form on $\Map(M,Z)$,
$$\varphi(X,Y)=\int_M \omega(X,Y)V$$
This is zero if $M$ is Lagrangian. Hence the $n$-dimensional orbits of
$T^n$ are Lagrangian themselves.
\vskip .25cm
If $X_1,\dots,X_n$ are the generating vector fields of $T^n=S^1\times
\dots\times S^1$, then since there are no complex directions in a
Lagrangian subspace, the vector fields $X_1,\dots,X_n$ together with
$IX_1,\dots,IX_n$ form a trivialization of the tangent bundle of ${\mathcal
M}$. Furthermore, since the $X_i$ preserve the complex structure $I$,
$$Z_1=X_1-iIX_1,\, \dots\dots,\, Z_n=X_n-iIX_n$$
are $n$ commuting holomorphic vector fields. In other words, a subset
${\mathcal M}_0\subset {\mathcal M}$  is biholomorphically equivalent to a
neighbourhood of
$$(S^1)^n\subset ({\C^*})^n$$
Choose ${\mathcal M}_0$ to be fibred over a ball ${\mathcal B}_0\subset {\mathcal B}$.
Note in passing that the trivialization by the holomorphic vector fields
defines a trivialization of the canonical bundle of ${\mathcal M}$.
\vskip .25cm
Now  ${\mathcal M}_0$  is homotopy equivalent to an orbit of $T^n$ and the K\"ahler form (which we shall call $\omega$) 
restricts to zero on each orbit since they are Lagrangian. Thus the
cohomology class  of $\omega$ is trivial, so there exists $\theta\in
\Omega^1({\mathcal M}_0)$ such that
$$\omega=d\theta=d(\theta^{1,0}+ \theta^{0,1})$$
Since $\omega$ is of type $(1,1)$, we have $\bar \partial \theta^{0,1}=0$
and so since the higher Dolbeault cohomology of a product of open sets in
${\C}^*$ vanishes, we can find a function $f$ such that $\bar\partial
f=\theta^{0,1}$, hence
$$\omega=d(\partial \bar f+\bar \partial f)=2i\partial \bar \partial \phi$$
so that $\phi=Im f$ is  a K\"ahler potential.
\vskip .25cm
Averaging over the compact group $T^n$ we can find a $T^n$-invariant
potential, and thus, using local coordinates
$$(z_1,\dots,z_n)\mapsto (e^{z_1},\dots,e^{z_n})\in {\C^*}^n$$
we have
$$\phi(z_1,\dots, z_n)=\phi(z_1+\bar z_1, \dots, z_n+\bar z_n)$$
The K\"ahler potential is thus a function only of the real part $x_i$ of
the holomorphic coordinates $z_i$. These metrics were studied by Calabi in
\cite{Cal}. In particular, the Ricci tensor vanishes and $dz_1\wedge\dots
\wedge dz_n$ is a covariant  constant holomorphic $n$-form if and only if
the function $\phi$ satisfies a real Monge-Amp\`ere equation
$$\det \frac{\partial^2 \phi}{\partial x_i\partial x_j}=const$$
The metric can now be written relative to the coordinates $x_i,y_i$ as
\begin{equation}
g=\sum_{i,j} \frac{\partial^2 \phi}{\partial x_i\partial x_j}(dx_idx_j+dy_idy_j)
\label{metric}
\end{equation}
\vskip .5cm
Since $T^n$ acts isometrically on ${\mathcal M}$, we can define a quotient
metric on the base space ${\mathcal B}$, the moduli space of special Lagrangian
submanifolds $M$. Since $x_1,\dots,x_n$ are $T^n$-invariant they define
local coordinates on ${\mathcal B}$ and the quotient metric has the Hessian form
$$g=\sum_{i,j} \frac{\partial^2 \phi}{\partial x_i\partial x_j}dx_idx_j$$
\vskip .25cm
There is another way of seeing this matrix of coefficients for the metric.
Note that in our local coordinates $z_i=x_j+iy_j$, the orbit of $T^n$ over
the point in ${\mathcal B}$ with coordinates $(x_1,\dots,x_n)$ is  $e^{iy_1},\dots
,e^{iy_n}$.   Thus
$$\frac{1}{2\pi}dy_1,\,\dots,\, \frac{1}{2\pi}dy_n$$
are forms whose de Rham cohomology classes form an integral basis for the
first cohomology of each fibre of ${\mathcal M}_0\rightarrow {\mathcal B}_0$. Now
$Idx_j=dy_j$, and $I$ is orthogonal so the metric on each torus is the flat
metric
$$g=\sum_{i,j} \frac{\partial^2 \phi}{\partial x_i\partial x_j}dy_idy_j$$
 Using that description, we have a map from the base space ${\mathcal B}_0$
into the moduli space of flat real $n$-tori, which can be identified with
the space of $n\times n$ positive definite matrices $\gl(n,\R)^+/\so(n)$.
This map is
$$x\mapsto \frac{\partial^2 \phi}{\partial x_i\partial x_j}$$
Globally it is well-defined modulo the action of $\gl(n,\Z)$, since
$dy_1/2\pi,\dots,dy_n/2\pi$ define an integral basis of $H^1(M,\R)$.
\vskip .25cm

\section{ MIRROR SYMMETRY}

\subsection{The SYZ approach}

In \cite{SYZ}, Strominger, Yau and Zaslow gave a geometrical approach to
mirror symmetry which has generated much interest recently. The reader is
referred in particular to the papers of M.Gross  \cite{Gr1},\cite{Gr2}. The
setting is that of a Calabi-Yau manifold $Z$ (let's take three dimensions
for simplicity) and a special Lagrangian 3-torus $M\subset Z$. As we know,
from McLean's theorem this belongs to a family of deformations through
special Lagrangian tori, whose dimension is
$$\dim {\mathcal B}=b_1(T^3)=3$$

\vskip .25cm
Strominger, Yau and Zaslow make the assumption that there is a smooth map
$p:Z\rightarrow B$ to a compact 3-manifold such that the generic fibre is
one of these special Lagrangian tori. Over some discriminant locus in $B$,
the fibres are supposed to degenerate into singular spaces, like the
elliptic curves in a holomorphic  elliptic fibration of an algebraic
surface. It has to be said that, at the moment of writing, no example of
this phenomenon is known which incorporates all of the data: a Calabi-Yau
metric, a complete family of degenerating special Lagrangian tori, and a
smooth base space. But then, the only way we have of finding Calabi-Yau
metrics on compact threefolds is by the existence theorem of Yau, and it is
not easy to imagine what a concrete example might be. Let us assume however
that such things exist, then we can make some observations.
\vskip .25cm
The first is that the complement $B'\subset B$ of the discriminant locus
can be identified with a global version of our moduli space ${\mathcal B}$ of
special Lagrangian deformations of a fixed torus. This is because both $B$ and $\mathcal B$ are three-dimensional so every sufficently small deformation in McLean's family is a fibre of $p$. We have one more property
from the fibration however: the normal vector fields which define the
deformations of a fibre $M$ over a point $x\in B'$ are just the pull-backs
of a basis
$$\frac{\partial}{\partial x_1},\, \frac{\partial}{\partial
x_2},\,\frac{\partial}{\partial x_3}$$
of the tangent space $(TB')_x$. In particular they are linearly independent
as sections of the normal bundle at each point of the fibre $M$. Hence the
three harmonic 1-forms
$$\iota(\frac{\partial}{\partial x_1})\omega,\,
\iota(\frac{\partial}{\partial x_2})\omega,\,\iota(\frac{\partial}{\partial
x_3})\omega$$
on the torus $M$ must be linearly independent at each point. This is not
necessarily true for an arbitrary metric on a torus but since it is true
for a flat torus it is at least true for metrics which are sufficiently
close to being flat.
\vskip .25cm
Given this special Lagrangian torus fibration, Strominger, Yau and Zaslow
propose that the mirror partner $\check Z$ of the Calabi-Yau manifold $Z$
should be the moduli space of pairs consisting of a special Lagrangian
torus and a flat line bundle over it. Clearly this has the same flavour as
the moduli space ${\mathcal M}$ we have described: pairs consisting of a
special Lagrangian torus and a degree one gerbe. To put everything in its
proper place, however, we need to study further properties of gerbes.

\subsection{Linear equivalence of points}

When 19th century mathematicians were studying divisors on Riemann
surfaces, they didn't use the language of line bundles, but nevertheless
established most of the basic theorems in the subject. Nowadays we would
say that given a set of points $p_1,\dots,p_n$  on a Riemann surface $X$,
the  divisor $D=p_1+\dots +p_n$ defines a line bundle $L_D$ with a section
which vanishes on $D$. Given another divisor $E=q_1+\dots +q_n$, then $D$
is linearly equivalent to $E$ if the two line bundles $L_D$ and $L_E$ are
holomorphically equivalent: $L_D\cong L_E\cong L$. If that is true then $D$
and $E$ are the zeros of two holomorphic sections $s_D,s_E$ of the same
line bundle $L$, so that $s_D/s_E$ is a meromorphic function with zeros on
$D$ and poles on $E$. And this is what our forebears used instead of line
bundles:  holomorphic functions on $X\backslash\{D\cup E\}$.
\vskip .25cm
Suppose we now have a set of points $p_1,\dots,p_n$ in a compact Riemannian
3-manifold $X$. We wish to introduce an analogous equivalence relation, but
what should we expect? The Riemann surface example can be thought of this
way: $s_D$ defines a trivialization of the line bundle $L$ outside $D$ and
$s_E$ outside $E$. The difference of the two trivializations is a function,
which for linear equivalence must be a holomorphic function. If $D$ and $E$
are now of codimension three, we should expect linear equivalence to be
framed in the language of differential equations on {\it line bundles}
defined in the complement of $D\cup E$ -- moving one step up in the
hierarchy. There is such a notion, which we shall see next.
\vskip .25cm
Recall from Lecture 1 that we defined a connection on the gerbe ${\mathcal
G}_p$ corresponding to a point $p\in X^3$ by using the Hodge theory of
currents. Its curvature was $V$ where the volume form was $V$ and  $X$ had
volume $2\pi$. So given two points, the gerbe ${\mathcal G}_p{\mathcal G}^{-1}_q$
has a connection whose curvature is zero. We shall say that {\it $p$ and
$q$ are linearly equivalent if the holonomy in $H^2(X,\R/\Z)$ of this flat
connection is trivial}. Similarly given two sets of points $p_1,\dots, p_n$
and $q_1,\dots,q_n$ we say that $p_1+\dots +p_n$ is linearly equivalent to
$q_1+\dots +q_n$ if the holonomy of ${\mathcal G}_{p_1}{\mathcal G}^{-1}_{q_1}\dots
{\mathcal G}_{p_n}{\mathcal G}^{-1}_{q_n}$ is trivial.
\vskip .25cm
Now the connection on ${\mathcal G}_p$ was obtained by taking an open  covering
of $X$ by a coordinate neighbourhood $U_p$ of $p$ and the complement of $p$
and putting $F_0=d^*H_p$ on $X\backslash \{p\}$ where
$$dd^*H_p=V-2\pi\delta_p$$
and choosing $F_1=d^*H_1$ with $dd^*H_1=V$ on $U_p$. Thus the
 connection on the gerbe ${\mathcal G}_p{\mathcal G}^{-1}_q$ is given by taking a
covering of $X$ by coordinate neighbourhoods $U_p$,$U_q$ of $p$ and $q$ and
$U=X\backslash\{p,q\}$ and then taking $F=d^*(H_p-H_q)$ on $U$,
$F_1=d^*(H_1-H_q)$ on $U_p$, $F_2=d^*(H_p-H_2)$ on $U_q$. These are all
smooth {\it closed} 2-forms on their respective open sets, and the
cohomology class of $F/2\pi$ restricted to $H^2(U_p\backslash
\{p\},\R)\cong \R$ or $H^2(U_q\backslash \{q\},\R)$ is integral.

Now it is easy to see that when we have a trivialization of a flat gerbe
(and a splitting $E\cong T^*\oplus 1$ as in Section 1.6) so that the gerbe
connection is represented by a  closed 2-form $F$, the holonomy, which we
described in Lecture 1 in \v Cech terms, is also given in de Rham terms by
the image of the  cohomology class
$$[F/2\pi]\in H^2(X,\R) \rightarrow H^2(X,\R/\Z)$$
This is directly analogous to a flat connection on a line bundle being
defined by a global closed 1-form $A$, whose periods modulo $2\pi \Z$
define the holonomy. But now in the Mayer-Vietoris sequence for $X=U\cup
(U_p\cup U_q)$ we have
$$0\rightarrow H^2(X, \R/\Z)\rightarrow H^2(U,\R/\Z)\rightarrow
H^2(U_p\backslash \{p\},\R/\Z)\oplus H^2(U_q\backslash \{q\},\R/\Z)$$
and so since the restriction of $[F/2\pi]\in H^2(U,\R)$  is an integral
class in  the cohomology groups $H^2(U_p\backslash \{p\},\R)$ and $H^2(U_q\backslash \{q\},\R)$,
the holonomy of the flat gerbe connection is uniquely determined by the
cohomology class of $F$ in $U$.
\vskip .25cm
It follows from the definition of linear equivalence that $p$ is linearly
equivalent to $q$ if and only if the closed form $F/2\pi$ defines an
integral class in $H^2(U,\R)$.
\vskip .25cm
We can interpret this result in the language of connections on line
bundles, as we hoped, since if $F/2\pi$ has integral periods, then $F$ is
the curvature of a connection $A$ on a line bundle defined over $U$. But
$F=d^*H$ where
$$dd^*H=2\pi(\delta_q-\delta_p)$$
so that $\ast H=\phi$ is a harmonic function with a singularity of the form
$1/2r$ at $q$ and $-1/2r$ at $p$. We have therefore a solution of the Dirac
monopole equation (or $S^1$ Bogomolny equations) on $X^3$:
$$F_A=\ast d_A\phi$$
This represents a charge 1 monopole located at $q$ and a charge $-1$
monopole at $p$. These are the natural gauge-theoretic equations for line
bundles in 3 dimensions.
\vskip .25cm
This discussion of linear equivalence using gerbe connections is
satisfactory in the sense that it sits nicely above line bundles and
divisors in the hierarchy, but it is not a practical way to determine if
sets of points are linearly equivalent. This is provided by the following
theorem:

\begin{theorem} (see \cite{Chat}) The sets of points $p_1,\dots, p_n$ and $q_1,\dots,q_n$ in a
compact Riemannian 3-manifold $X$ are linearly equivalent if and only if
$$\sum_i\int^{q_i}_{p_i} \theta \in \Z$$
for each harmonic 1-form $\theta$ on $X$ with integral cohomology class.
\end{theorem}
The contour integral in the theorem is taken by using any smooth path
$\gamma_i$ from $p_i$ to $q_i$. Since the integrals of $\theta$ over closed
paths are integral, the condition is independent of the path.
\vskip .25cm
\begin{proof} To address the theorem, let us for simplicity take just two
points $p$ and $q$. Then these are linearly equivalent if and only if the
2-form $F/2\pi$ has integral periods in $U=X\backslash\{p,q\}$, where
$F=d^*H$ and
$$dd^*H=2\pi (\delta_q-\delta_p)$$
If $\varphi$ is a closed 1-form with compact support in $U$ and integral
cohomology class in the compactly supported cohomology of $U$, then this is
equivalent to the condition
\begin{equation}
\int_U \varphi \wedge F\in 2\pi \Z
\label{periods}
\end{equation}
for all such $\varphi$.
Now let $\gamma$ be a path from $p$ to $q$. As a current we  define
$$\langle \gamma, \theta \rangle =\int_{\gamma} \theta$$
but as a current $\gamma$ also has a Hodge decomposition \cite{dR}
\begin{equation}
\gamma=A+dB+d^*C
\label{hodge}
\end{equation}
where $A$ is the harmonic part, which by elliptic regularity is a smooth
globally defined harmonic 2-form. Now by Stokes' theorem
$$\langle d\gamma, f \rangle=\langle \gamma, df \rangle=\int_{\gamma}
df=f(q)-f(p)$$
so that
$$d\gamma=\delta_q-\delta_p$$
Hence from (\ref{hodge})
$$dd^*C=\delta_q-\delta_p=\frac{1}{2\pi}dd^*H$$
and so $2\pi d^*C=d^*H=F$.
\vskip .25cm
Now since the boundary of $\gamma$ is $q-p$, $\gamma$ is an integral
relative cycle in $U=X\backslash\{p,q\}$, and so if $\varphi$ is a
compactly supported 1-form with integral cohomology class,
$$\int_{\gamma} \varphi \in \Z$$
Now since $d\varphi=0$, $\langle dB,\varphi\rangle=0$ so from (\ref{hodge})
$$\langle \gamma, \varphi\rangle=\langle A,
\varphi\rangle+\frac{1}{2\pi}\langle F,\varphi\rangle$$
and hence since $\langle \gamma, \varphi\rangle\in \Z$, $\langle
F,\varphi\rangle\in 2\pi \Z$ if and only if $\langle A, \varphi\rangle\in
\Z$.
\vskip .25cm
Now let $\theta$ be a harmonic 1-form on $X$, then since $d\theta=d^*\theta=0$,
$$
\langle \gamma,\theta\rangle =\langle A+dB+d^*C, \theta \rangle=\langle
A,\theta\rangle
$$
But $\theta$ is cohomologous to a closed form $\varphi$  with compact
support in $U$: all we have to do is solve $df=\theta$ in the open sets
$U_p$ and $U_q$, extend to a global function using a partition of unity and
write $\varphi=\theta-df$. It follows, since $A$ is closed, that
$$\langle A,\theta\rangle=\langle A,\varphi\rangle$$

Thus finally, $p$ and $q$ are linearly equivalent if and only if $\langle
A,\theta\rangle$ is an integer for all harmonic 1-forms with integral
cohomology class.
\end{proof}
\vskip .5cm
Another way to formulate Theorem 3.1 is to define the Abel-Jacobi map
$$u:X\rightarrow H^2(X,\R/Z)$$
by choosing a base point $p$ and defining
$$\langle u(x), [\theta]\rangle=\int_p^x\theta$$
for $[\theta]\in H^1(X,\R)$ the de Rham cohomology class of the harmonic
1-form $\theta$. Then $p_1+\dots+p_n$ is linearly equivalent to
$q_1+\dots+q_n$ if and only if
$$\sum_i u(p_i)=\sum_i u(q_i)$$
\subsection{Degree one gerbes}
Recall our K\"ahler quotient construction in Lecture 2 of a natural
K\"ahler metric on the moduli space ${\mathcal M}$ of special Lagrangian
submanifolds and degree one gerbes. As we have already noted, in a special
Lagrangian torus fibration this is 3-dimensional, just like the Calabi-Yau
manifold $Z$. There is a close relationship, as we shall see now.

Let $Z'\subset Z$ be the open set of $Z$ which is fibred by nonsingular
special Lagrangian tori. Then  each point $z\in Z'$ lies on a unique
special Lagrangian torus $M_z$. Furthermore $z\in M_z$ defines a degree one
gerbe ${\mathcal G}_z$ on $M_z$ where we take the normalized volume form for
its curvature. We therefore obtain a natural map:
$$J: Z'\rightarrow {\mathcal M}$$
defined by $J(z)=(M_z, {\mathcal G}_z)$.

\begin{theorem} Let $Z$ be a 3-dimensional Calabi-Yau manifold with a special
Lagrangian torus fibration. Then the map $J: Z'\rightarrow {\mathcal M}$ is a
diffeomorphism.
\end{theorem}

\begin{proof} The map $J$ is clearly smooth and commutes with the
projections onto ${\mathcal B}$, the moduli space of special Lagrangian tori,
so it remains to prove that the map on the fibres is a diffeomorphism.

Consider a fibre $M\cong T^3\subset Z'$. The corresponding fibre on ${\mathcal
M}$ is the moduli space of degree one gerbes on $M$, also a 3-torus.
If we choose a basepoint $p\in M$, and make $[{\mathcal G}_p]$ an origin in the
moduli space of gerbes, then the map $J$ can be interpreted as the map
$$j:M\rightarrow H^2(M,\R/\Z)$$
where $j(q)$ is the holonomy of the flat gerbe ${\mathcal G}_q{\mathcal G}^{-1}_p$.
We can easily see from the proof of Theorem 3.1 that this is the
Abel-Jacobi map
$u$ (in particular it is clear that $j(p)=j(q)$ if and only if $p$ and $q$
are linearly equivalent).
\vskip .25cm
Let us consider the kernel of the derivative of $j=u$ at $q\in M$. By the
definition of $u$, this consists of the tangent vectors $X\in T_qM$ such
that
$$\iota(X)\theta\vert_q=0$$
for all harmonic 1-forms $\theta$ on $M$. But as we noted, for a fibration
the harmonic 1-forms are linearly independent at each point so the kernel
is always zero and $j$ is a local diffeomorphism. Since $M$ is compact it
is a covering map.
\vskip .25cm
 \noindent But now the Abel-Jacobi map is defined for all metrics, and for
the flat metric on the torus, $u:T^3\rightarrow H^2(T^3,\R/\Z)$ is a
diffeomorphism. In particular the induced map on $\pi_1(T^3)=H_1(T^3,\Z)$
is an isomorphism. Any metric can be continuously connected to the flat
metric hence the induced map:
$$u_*:H_1(T^3,\Z)\rightarrow H_1(H^2(T^3,\R/\Z),\Z)$$
 is an isomorphism for all metrics. It follows that the covering map is a
diffeomorphism and the theorem is proved.
\end{proof}
\vskip .25cm
As a consequence of this result,  we  now see $Z'$ naturally (but not
isometrically or symplectically) identified with the  moduli space ${\mathcal
M}$ of special Lagrangian tori and degree one gerbes.
We have given in (\ref{metric}) a local analytical form for the natural
metric on ${\mathcal M}$
$$g=\sum_{i,j} \frac{\partial^2 \phi}{\partial x_i\partial
x_j}(dx_idx_j+dy_idy_j)$$
 obtained by a quotient construction, but it will be useful to see some of
its features in a more geometrical fashion. Note that if $T_F$ is the
tangent bundle along the fibres, the orbits of $T^3$, then $IT_F$ is
transversal to the fibres and $T^3-$invariant. It therefore defines the
horizontal subspaces of a {\it connection} on ${\mathcal M}$ considered as a
principal $T^3$ bundle over ${\mathcal B}$. This connection is moreover flat,
since the horizontal subbundle is spanned by the commuting vector fields
$IX_1,IX_2,IX_3$. From the form of the metric, in the local coordinates
$x_i,y_i$ the horizontal subspaces integrate to the submanifolds
$(y_1,y_2,y_3)=const$.
\vskip .25cm
There is also a flat connection arising from the point of view of gerbes.
Consider the 3-form $\Omega_1$ on $Z$ and restrict to $Z_0=p^{-1}({\mathcal
B}_0)$  By assumption the closed form $\Omega_1/2\pi$ restricts to an
integral class on each fibre, but $Z_0$ is homotopy equivalent to a fibre,
so the integrality condition holds on  $Z_0$ itself. We can therefore find
a gerbe ${\mathcal H}$ with connection and curvature $\Omega_1$ on $Z_0$ which
restricts to a degree one gerbe on each fibre of $Z_0$. Any two such
choices differ by a flat gerbe whose holonomy in $H^2(Z_0, \R/\Z)$ is
determined (again by homotopy invariance) by its value on any fibre. Thus
choosing such gerbes with connection on $Z_0$ gives a family of sections of
${\mathcal M}_0\rightarrow {\mathcal B}_0$ and their tangent spaces generate a flat
connection on the whole of ${\mathcal M}$. This is a Gauss-Manin connection in
the context of gerbes. We shall show that it coincides with the connection
$IT_F$.
\vskip .25cm
A horizontal space of the Gauss-Manin connection integrates to the
equivalence classes  of maps $f:M\rightarrow Z$ and degree one gerbes
${\mathcal G}$ such that $$ {\mathcal G}\cong f^*{\mathcal H}$$
Fix three 2-dimensional tori $T_1,T_2,T_3$ in $M$ whose homology classes
form generators of  $H_2(M,\Z)$. For each gerbe ${\mathcal G}$, restrict to $T_i$. Since
the subtori are 2-dimensional the gerbe is flat on each of them, so take
the holonomy  on each to  define a map
$$H:{\mathcal M}\rightarrow (\R/\Z)^3$$
Tensoring with a flat gerbe on $M$, we see that the map $H$ is equivariant
with respect to the action of $T^3$ and surjective, so a fibre is defined by
$(y_1,y_2,y_3)=(c_1,c_2,c_3)$. But by definition, a flat leaf of the
Gauss-Manin connection is the subset for which  $ {\mathcal G}=f^*{\mathcal H}$ has
fixed holonomy around the 2-tori, i.e. is given by $(y_1,y_2,y_3)=constant$
and so we see that the two flat connections coincide: one defined by the
complex structure $I$, and the other by homotopy invariance for the
parametrized objects.
\vskip .25cm
The result we have here -- an identification of the original Calabi-Yau
with a moduli space of special Lagrangian tori and gerbes -- is not part of
the original Strominger-Yau-Zaslow approach but it puts the SYZ mirror into
the same context as the original Calabi-Yau. There is a difference however.
Whereas, our moduli space ${\mathcal M}$
 consists of a bundle over ${\mathcal B}$ whose fibres are flat tori isomorphic
to  $H^2(M,\R/\Z)$, the SYZ mirror is a similar bundle whose fibres are
isomorphic to the moduli space of flat connections on $M$ -- the {\it dual
torus} $H^1(M,\R/\Z)$. Before we continue to relate the two, we need to
consider linear equivalence in a different setting, for codimension three
submanifolds of the Calabi-Yau.

\subsection{Linear equivalence of special Lagrangian submanifolds}

We saw right at the beginning that topologically we can associate a gerbe
to an oriented submanifold of codimension 3. Since $\Omega_1$ restricts to
a volume form on a special Lagrangian submanifold, then special Lagrangians
in 3-dimensional Calabi-Yau manifolds are codimension three and oriented.
We need to consider the differential geometry of the gerbes they generate.

The situation for general codimension three submanifolds is very similar to
that of points in a 3-manifold.
Given $M^{n-3}\subset X^n$ and $G\in \Omega^3(X)$  a closed 3-form
cohomologous to $2\pi M$ in the de Rham cohomology of currents, we can
solve the equation for currents
$$\Delta H=G-2\pi M$$
In this case $H$ is a singular form of degree $(n-3)$, but because
$dG=dM=0$, $\Delta dH=0$ and $dH$ is harmonic. But it is  also exact and
hence zero, so as before we have
$$dd^*H=G-2\pi M$$
Now on coordinate neighbourhoods $U_\alpha$ of $M$ we can solve
$dd^*H_\alpha=G$ with $dH_\alpha=0$ (this is basic local potential theory).
We cover $X$ with these open sets and $U_0=X\backslash N(M)$. Then on
$U_{0}\cap U_{\alpha}$, $d(d^*H-d^*H_\alpha)=0$ and has  cohomology class
in $2\pi\Z$
so is the curvature $F$ of  a connection $\nabla_{\alpha 0}$. This way we
construct the data necessary to give a gerbe connection: line bundle
connections
$\nabla_{\alpha\beta}, \nabla_{\alpha 0}$
with curvature
$$ F_{\alpha 0}=d^*H-d^*H_\alpha,\quad F_{\alpha\beta}=d^*H_\beta-d^*H_\alpha$$
and 2-forms
$$F_0=d^*H,\quad F_\alpha=d^*H_\alpha$$
such that $F_{\beta}-F_{\alpha}=F_{\alpha\beta}$. The
 curvature of the connection is $G$ (since $dd^*H=G=dd^*H_\alpha)$.
 \vskip .25cm
 Since a codimension three submanifold $M$ defines a gerbe ${\mathcal G}_M$
with connection, we can ask for linear equivalence in the same sense as for
points: if $M^{n-3}$ and $N^{n-3}$ are homologous, we say that they are
linearly equivalent if the flat gerbe ${\mathcal G}^{-1}_M{\mathcal G}_N$ has
trivial holonomy. And again, using the language of currents, $M$ and $N$
are linearly equivalent if and only if
$$\int_\Gamma\theta \in   \Z$$
for each harmonic $(n-2)$-form $\theta$ on $X$ with integral cohomology
class, where $\Gamma$ is a smooth chain with $\partial \Gamma=M-N$, for
example an $(n-2)$ manifold with boundary $M-N$.
\vskip .25cm
Let us focus now on special Lagrangian submanifolds of a 3-dimensional
Calabi-Yau manifold. Each such manifold $M^3\subset Z$ has, as we have
seen, a unique family of deformations as special Lagrangian submanifolds.
An obvious question is: are these deformations all linearly equivalent to
$M$? There is a straightforward answer:

\begin{theorem} A special Lagrangian submanifold in a compact Calabi-Yau
3-manifold is linearly equivalent to any of its deformations if and only if
the restriction map  $H^2(Z, \R)\rightarrow H^2(M,\R)$ is zero.
\end{theorem}

\begin{proof} Let $\gamma$ be a curve in the moduli space ${\mathcal B}$ of
deformations from the point $p\in {\mathcal B}$ representing $M$ and an
arbitrary point $q\in {\mathcal B}$. Then there is a corresponding smooth map
$$F:M\times [0,1]\rightarrow Z$$
such that $F_t(M)$ is the deformation of $M=F_0(M)$. Now $M$ is linearly
equivalent to $F_t(M)$ if and only if
$$\int_{M\times [0,t]}F^*\theta \in \Z$$
for each harmonic $4$-form $\theta \in \Omega^4(Z)$ with integral
cohomology class. Differentiating with respect to $t$ at $t=0$ this implies
that
\begin{equation}
\int_M\iota(X)\theta=0
\label{int}
\end{equation}
where $X$ is the section of the normal bundle defining an infinitesimal
variation of special Lagrangian submanifolds. This must hold for all
harmonic $\theta$ with integral cohomology class and so by linearity for
all $\theta$.
\vskip .25cm
Now by Hodge theory the map $L:\Omega^2(Z)\rightarrow \Omega^4(Z)$ defined
by $L\alpha=\omega\wedge \alpha$ is an isomorphism from harmonic 2-forms to
harmonic 4-forms, so we can write
$$\theta=\omega \wedge \alpha$$
for a harmonic 2-form $\alpha$. Now
$$\iota(X)\theta=(\iota(X)\omega) \wedge \alpha+\omega\wedge (\iota(X)\alpha)$$
but $M$ is Lagrangian, so that $\omega$ restricts to zero on $M$, thus
(\ref{int}) becomes
$$\int_M(\iota(X)\omega) \wedge \alpha=0$$
But now we know that $\iota(X)\omega$ restricted to $M$ is a harmonic form,
and moreover any harmonic form is the derivative of a variation through
special Lagrangian submanifolds. Thus if all deformations are linearly
equivalent to $M$,
$$ \int_M\theta\wedge \alpha=0$$
for all harmonic 2-forms $\alpha$ on $Z$ and all harmonic 1-forms $\theta$
on $M$. Since the harmonic 1-forms run through each cohomology class this
means that the cohomology class $[\alpha]$ of $\alpha$ on $M$ vanishes,
which is the theorem. The converse follows by differentiating the integral
$$\int_{M\times [0,t]}F^*\theta $$
with respect to $t$.
\end{proof}
\vskip .25cm
If two codimension 3 submanifolds are linearly equivalent on a
3-dimensional Calabi-Yau, then as before we have a line bundle with
connection and curvature $F$ on the complement of the two submanifolds, and
such that the 4-form $\ast F$ is closed. If the submanifolds are special
Lagrangian we can say more.
\vskip .25cm
First consider the Lagrangian condition that $\omega$ vanishes on $M$. This
means that
$$\int_M \omega\wedge \theta=0$$
for all 1-forms $\theta$, so as a current $\omega\wedge M=0$. So if
$\Delta H=M-N$ and both $M$ and $N$ are special Lagrangian,
$$\Delta (H\wedge \omega)=0$$
Thus $H\wedge \omega$ is a smooth (by elliptic regularity) harmonic 5-form.
In fact on a compact Calabi-Yau such forms exist only if they are covariant
constant, so if we assume that the metric is irreducible then $H\wedge
\omega=0$, so that $H$ is a primitive 3-form.
\vskip .25cm
We also have the special Lagrangian condition that $\Omega_2$ restricts to
zero on $M$ and $N$. This means similarly that
$$\Delta (H\wedge \Omega_2)=0$$
so that $(H\wedge \Omega_2)$ is a constant multiple of the volume form of
$Z$. In fact, if we decompose the current $H$ as
$$H=\phi\Omega_1+b\Omega_2+h$$
where $h$ is a current of type $(2,1)+(1,2)$ then the above observations
say that $h$ is  primitive and $b$ is constant. Now the curvature $F$ of
our connection is given by $d^*H$. It follows from this that $F$ satisfies
the conditions
$$\Lambda F^{1,1}=0,\quad F^{2,0}=\ast_3 \partial \phi$$
where $\ast_3:\Omega^{1,0}(Z)\rightarrow \Omega^{2,0}(Z)$ is the complex
Hodge star operator on Calabi-Yau manifolds introduced by Donaldson and
Thomas in \cite{DT}. In fact these equations are the dimensional reductions
to 3 dimensions of their (abelian) $\su(4)$ instanton equations on
Calabi-Yau 4-manifolds.

\subsection{ SYZ with a  B-field}

Now let us return to the Strominger-Yau-Zaslow construction of the mirror
$\check Z$ of a 3-dimensional Calabi-Yau manifold $Z$ with a special
Lagrangian torus fibration. Their mirror is the moduli space of pairs
$(M,L)$ where $M$ is a special Lagrangian torus in the fibration and $L$ a
flat line bundle over $M$. If we compare this to our moduli space ${\mathcal
M}$ which models $Z$, then we notice a difference. Both are torus
fibrations, but the torus of degree 1 gerbes has no distinguished origin,
whereas the class of the trivial flat line bundle gives the mirror a
distinguished section. We can balance the picture by introducing another
object, in fact a flat gerbe, which physicists do in any case for other
reasons. This is the {\it B-field}.
\vskip .25cm
Very briefly, mirror symmetry is supposed to exhibit certain features, one
of which is a duality between the Dolbeault spaces $H^1(Z, T)$ and
$H^1(\check Z, T^*)$. The first defines infinitesimal deformations of the
complex structure, the second deformations of the K\"ahler class
$[\omega]$. This is fine except that  $[\omega]$ is real, and indeed the
space $H^1(\check Z, T^*)$ has a real structure from complex conjugation on
$(1,1)$ forms whereas $H^1(Z, T)$ is complex. A twin field to the
cohomology class of the K\"ahler metric is introduced to restore the
balance, and this is the  B-field. What exactly it is geometrically may not
be  clear in the literature, but the fact that it carries a cohomology
class in $H^2(Z,\R/\Z)$ is important. Given the general drift of our
approach in these lectures, it is reasonable to suggest that the  B-field
is in fact a flat gerbe ${\bf B}$ and the class in $H^2(Z,\R/\Z)$ is its
holonomy. Let us make that assumption, then we can define an analogue of
the SYZ mirror where the torus fibration has no distinguished zero section.

\vskip .25cm
Let us make a further assumption that the special Lagrangian torus fibres
are linearly equivalent, then from Theorem 3.3, the restriction map for the second cohomology 
$H^2(Z,\R)\rightarrow H^2(M,\R)$ is zero. This means that the restriction
map $H^2(Z,\R/\Z)\rightarrow H^2(M,\R/\Z)$ is trivial and so  the flat
gerbe ${\bf B}$ has trivial holonomy on each torus fibre. We can therefore
take its flat trivializations, any two of which differ by a flat line
bundle. We thus define the SYZ mirror $\check Z$ of a Calabi-Yau with a
B-field as {\it the moduli space of pairs $(M,T)$ where $M$ is a special
Lagrangian torus and $T$ is a flat trivialization of the gerbe ${\bf B}$ on
$M$}. When ${\bf B}$ is the trivial flat gerbe, we obtain the original SYZ
mirror.
\vskip .25cm
We also have a Gauss-Manin connection for this definition of the mirror.
The inverse in $\check Z$ of a small ball in ${\mathcal B}$ is homotopy
equivalent to a fibre, and so the holonomy of ${\bf B}$ vanishes there too,
and any flat trivialization is determined by its restriction to a fibre. By
comparing with the metric on ${\mathcal M}\cong Z$, we use this flat connection
to split the tangent bundle into horizontal and vertical parts. Compare the
two:
$$TZ\cong H^2(M,\R)\oplus H^1(M,\R),\quad T{\check Z}\cong H^1(M,\R)\oplus
H^1(M,\R)$$
Here we see that, without any further information, the natural pairing on
cohomology defines the {\it symplectic form} $\omega$ on $Z$. By contrast
we have
$$T{\check Z}\cong H^1(M,\R)\otimes {\C}$$
so that we naturally get an {\it almost complex structure}. There is more
than that however. We use the splitting for $T{\check Z}$ to define a
metric on $\check Z$. On the horizontal space we put the metric from the
base ${\mathcal B}$ and on the fibre the metric of the torus $H^1(M,\R/\Z)$.
This is the dual torus to $H^2(M,\R/\Z)$ so the metric is
\begin{equation}
\check g= \sum_{i,j} g_{ij}dx_idx_j+g^{ij}d\eta_id\eta_j
\label{check}
\end{equation}
where
$$g_{ij}=\frac{\partial^2\phi}{\partial x_i\partial x_j}$$
and $g^{ij}$ is the inverse matrix of $g_{ij}$ representing the  metric on
the dual torus. This is not obviously a K\"ahler metric, but  following
\cite{H} we shall see that it is, in fact.
\vskip .25cm
Let $V$ be a real vector space and $V^*$ its dual, then we can define a
constant symplectic form on $V\times V^*$ by
$$\omega((a,\alpha),(b,\beta))=\langle a,\beta\rangle-\langle b,\alpha\rangle$$
We can also define an indefinite flat metric by
$$g((a,\alpha),(a,\alpha))=\frac{1}{2}\langle a, \alpha\rangle$$
If we think of $V\times V^*$ as the cotangent bundle of $V$, then $\omega$
is just the canonical symplectic form $\omega = \sum dx_i\wedge d\xi_i$ and
a general Lagrangian submanifold $Y\subset V\times V^*$  is the graph of
the derivative of a function:
$$\xi_i=\frac{\partial f}{\partial x_i}$$
The induced metric on $Y$ is
$$\sum_i dx_id\xi_i=\sum_{i,j} \frac{\partial^2 f}{\partial x_i\partial
x_j}dx_idx_j=\sum_{i,j}g_{ij}dx_idx_j$$
This is in precisely the Hessian form of the metric on ${\mathcal B}$.
\vskip .25cm
On the other hand, we can think of $V\times V^*$ as the cotangent bundle
$T^*V^*$, in which case  $Y$ can be written as
$$x_i=\frac{\partial \check f}{\partial \xi_i}$$
But then
$$\delta_{ij}=\sum_k \frac{\partial^2 \check f}{\partial \xi_i \partial
\xi_k}\frac{\partial \xi_k}{\partial x_j}=\sum_k \frac{\partial^2 \check
f}{\partial \xi_i \partial \xi_k}\frac{\partial^2 f}{\partial x_j \partial
x_k}$$
so that
$$g^{ij}=\frac{\partial^2 \check f}{\partial \xi_i \partial \xi_j}$$
and the same induced metric on $Y$ is
$$\sum_{i,j} \frac{\partial^2 f}{\partial x_i\partial
x_j}dx_idx_j=\sum_{i,j} \frac{\partial^2 \check f}{\partial \xi_i\partial
\xi_j}d\xi_id\xi_j$$
The transformation $f\mapsto \check f$ is the classical Legendre transform.
\vskip .25cm
\noindent Thus, if we use coordinates $\xi_j,\eta_j$ instead of
$x_j,\eta_j$ the metric (\ref{check}) becomes
$$\check g=\sum_{ij}\frac{\partial^2 \check \phi}{\partial \xi_i\partial
\xi_j}(d\xi_id\xi_j+d\eta_id\eta_j)$$
This is, as before, K\"ahler with potential $\check \phi$.
\vskip .25cm
Thus the SYZ approach coupled with the theory of gerbes, gives a symmetry
between metrics on $Z$ and $\check Z$. Neither of these metrics is part of
a metric on a true compact Calabi-Yau manifold, because of the high degree
of symmetry: there are no Killing fields on an irreducible compact
Riemannian manifold with zero Ricci tensor. The most one can hope for is
that this metric is the leading order term in a more general expression
which incorporates what are known as ``instanton corrections".

\bibliographystyle{amsalpha}

\begin{thebibliography}{11}
%
\bibitem{A}
M.F. Atiyah, \textit{Complex analytic connections in fibre bundles}, Trans. Am.
Math. Soc. \textbf{85} (1957), 181-207.
%
\bibitem{B}
R. L. Bryant,
\textit{Some examples of special Lagrangian tori}, preprint \textbf{math/9902076}
%
\bibitem{Br}
J.-L.  Brylinski, \textit {Characteristic classes and geometric quantization},
Progr. in Mathematics \textbf{107}, Birkh\"auser, Boston (1993)
%
\bibitem{Cal}
E. Calabi,
\textit{A construction of nonhomogeneous Einstein metrics},
 Proc. of Symp. in Pure Mathematics,
\textbf{27}   17-24,  AMS, Providence  (1975)
%
\bibitem{Chat}
D. S. Chatterjee, \textit{ On the construction of abelian gerbs}, PhD thesis
(Cambridge) (1998)
%
\bibitem{Don}
S. K. Donaldson,
\textit{Moment maps and diffeomorphisms},  Asian Journal of Mathematics, (to
appear)
%
\bibitem{DT}
S. K. Donaldson and R. P. Thomas,
\textit{Gauge theory in higher dimensions},
in \textit {The Geometric Universe: Science, Geometry and the work of Roger
Penrose}, S.A.Huggett et al (eds.), Oxford University Press, Oxford, 1998.
%
\bibitem{Gir}
J. Giraud, \textit {Cohomologie non-ab\'elienne}, Grundl. \textbf{179}, Springer
Verlag, Berlin  (1971)
%
\bibitem{Gr1}
M. Gross,
\textit{Special Lagrangian fibrations I: Topology},
preprint \textbf{alg-geom/9710006}
%
\bibitem{Gr2}
M. Gross,
\textit{Special Lagrangian fibrations II: Geometry},
preprint \textbf {math/9809072}
%
\bibitem{H}
N. J. Hitchin,
 \textit{The moduli space of special Lagrangian submanifolds},
 Dedicated to Ennio De Giorgi.  Ann. Scuola
Norm. Sup. Pisa Cl. Sci.  \textbf {25} (1997), 503--515
%
\bibitem{Mac}
R. C. McLean,
 \textit{Deformations of calibrated submanifolds},  Comm. Anal. Geom. \textbf {6}
(1998), 705--747
\bibitem{Mo}
J. Moser.
\textit{On the volume elements on a manifold},
       Trans. Amer. Math. Soc. \textbf {120} (1965) 286--294
%
\bibitem{dR}
 G. de Rham, {\it Differentiable manifolds}, Grundl. \textbf {266},
Springer Verlag, Berlin  (1984)
%
\bibitem{Sten}
 M. B. Stenzel,
  \textit{Ricci-flat metrics on the complexification of a compact rank
one symmetric space},  Manuscripta Math. \textbf {80} (1993), 151--163.
%
\bibitem{SYZ}
A. Strominger, S-T. Yau, E. Zaslow,
\textit{Mirror Symmetry is T-duality},  Nuclear Phys. \textbf {B 479} (1996),  243--259.
%
\bibitem{Wood}
N. M. J. Woodhouse, \textit {Geometric quantization}, Oxford University Press,
Oxford (1980).

\end{thebibliography}

\end{document}